\documentclass[journal, onecolumn]{IEEEtran}

\usepackage{graphicx}          
\usepackage{amsmath} 
\usepackage[english]{babel}
\usepackage{amsfonts}
\usepackage{amstext}
\usepackage{caption}
\usepackage{subcaption}
\usepackage{tabularx}
\usepackage{savesym}
\savesymbol{AND}
\usepackage{algorithm}
\usepackage{algorithmic}
\usepackage{color}

\usepackage{exStyle}

\newtheorem{theorem}{Theorem}[section]
\newtheorem{lemma}[theorem]{Lemma}
\newtheorem{proposition}[theorem]{Proposition} 
\newtheorem{corollary}[theorem]{Corollary}

\newtheorem{assumption}[theorem]{Assumption}
\newtheorem{remark}[theorem]{Remark}
\newtheorem{problem}[theorem]{Problem}

\begin{document}
%

\title{Inexact Alternating Minimization Algorithm for Distributed Optimization with an Application to Distributed MPC}

\author{Ye Pu, Colin N. Jones and Melanie N. Zeilinger %
\thanks{Y. Pu and C.N. Jones are with the Automatic Control Lab, \'Ecole Polytechnique F\'ed\'erale de Lausanne, EPFL-STI-IGM-LA
Station 9 CH-1015 Lausanne, Switzerland, e-mail: {\tt \small \{y.pu,colin.jones\}@epfl.ch}. }%
\thanks{M.N. Zeilinger is with the Empirical Inference Department, Max Planck Institute for Intelligent Systems, 72076 T\"ubingen, Germany, e-mail: {\tt\small melanie.zeilinger@tuebingen.mpg.de}. }%
\thanks{This work has received funding from the European Research Council under the European Union's Seventh Framework Programme (FP/2007-2013)/ ERC Grant Agreement n. 307608. The research of M. N. Zeilinger has received funding from the EU FP7 under grant agreement no. PIOF-GA-2011-301436-``COGENT''.}%
}




\maketitle

%
\IEEEpeerreviewmaketitle

\begin{abstract}

In this paper, we propose the inexact alternating minimization algorithm (inexact AMA), which allows inexact iterations in the algorithm, and its accelerated variant, called the inexact fast alternating minimization algorithm (inexact FAMA). We show that inexact AMA and inexact FAMA are equivalent to the inexact  proximal-gradient method and its accelerated variant applied to the dual problem. Based on this equivalence, we derive complexity upper-bounds on the number of iterations for the inexact algorithms. We apply inexact AMA and inexact FAMA to distributed optimization problems, with an emphasis on distributed MPC applications, and show the convergence properties for this special case. By employing the complexity upper-bounds on the number of iterations, we provide sufficient conditions on the inexact iterations for the convergence of the algorithms. We further study the special case of quadratic local objectives in the distributed optimization problems, which is a standard form in distributed MPC. For this special case, we allow local computational errors at each iteration. By exploiting a warm-starting strategy and the sufficient conditions on the errors for convergence, we propose an approach to certify the number of iterations for solving local problems, which guarantees that the local computational errors satisfy the sufficient conditions and the inexact distributed optimization algorithm converges to the optimal solution.

\end{abstract}

%

\section{Introduction}\label{sec:Introduction}

First-order optimization methods, see e.g. \cite{nesterov_method_1983}, \cite{beck_fISTA_2009} and \cite{goldstein_fast_2012}, play a central role in large-scale convex optimization, since they offer simple iteration schemes that only require information of the function value and the gradient, and have shown good performance for solving large problems with moderate accuracy requirements in many fields, e.g. optimal control \cite{richter_computational_2012}, signal processing \cite{combettes_proximal_2011} and machine learning \cite{boyd_distributed_2011}. In this paper, we will study a sub-group of first-order methods, called splitting methods, and apply them to distributed optimization problems. Splitting methods, which are also known as alternating direction methods, are a powerful tool for general mathematical programming and optimization. A variety of different spitting methods exist, requiring different assumptions on the problem setup, while exhibiting different properties, see e.g. \cite{goldstein_fast_2012} and \cite{combettes_proximal_2011} for an overview. The main concept is to split a complex convex minimization problem into simple and small sub-problems and solve them in an alternating manner. For a problem with multiple objectives, the main strategy is not to compute the descent direction of the sum of several objectives, but to take a combination of the descent directions of each objective. 

The property of minimizing the objectives in an alternating way provides an efficient technique for solving distributed optimization problems, which arise in many engineering fields \cite{boyd_distributed_2011}. By considering the local cost functions, as well as local constraints, as the multiple objectives of a distributed optimization problem, splitting methods allow us to split a global constrained optimization problem into sub-problems according to the structure of the network, and solve them in a distributed manner. The advantages of using distributed optimization algorithms include the following three points: in contrast to centralized methods, they do not require global, but only local communication, i.e., neighbour-to-neighbour communication; secondly, they parallelize the computational tasks and split the global problem into small sub-problems, which reduces the required computational power for each sub-system; thirdly, distributed optimization algorithms preserve the privacy of each-subsystem in the sense that each sub-system computes an optimal solution without sharing its local cost function and local constraint with all the entities in the network.

In this paper, we consider a distributed Model Predictive Control problem as the application for the distributed optimization to demonstrate the proposed algorithms, as well as the theoretical findings. Model Predictive Control (MPC) is a control technique that optimizes the control input over a finite time-horizon in the future and allows for constraints on the states and control inputs to be integrated into the controller design. However, for networked systems, implementing an MPC controller becomes challenging, since solving an MPC problem in a centralized way requires full communication to collect information from each sub-system, and the computational power to solve the global problem in one central entity. Distributed model predictive control \cite{Scattolini2009} is a promising tool to overcome the limiting computational complexity and communication requirements associated with centralized control of large-scale networked systems. The research on distributed MPC has mainly focused on the impact of distributed optimization on system properties such as stability and feasibility, and the development of efficient distributed optimization algorithms. 

\color{black}
However, a key challenge in practice is that  distributed optimization algorithms, see e.g. \cite{bertsekas_parallel_2003}, \cite{boyd_distributed_2011} and \cite{giselsson_accelerated_2013}, may suffer from inexact local solutions and unreliable communications. The resulting inexact updates in the distributed optimization algorithms affect the convergence properties, and can even cause divergence of the algorithm. 

\color{black}

In this work, we study inexact splitting methods and aim at answering the questions of how these errors affect the algorithms and under which conditions convergence can still be guaranteed. Seminal work on inexact optimization algorithms includes \cite{necoara_rate_2013}, \cite{dinh_fast_2012}, \cite{lin_universal_2015} and \cite{schmidt_convergence_2011}. In \cite{necoara_rate_2013}, the authors studied the convergence rates of inexact dual first-order methods. In \cite{dinh_fast_2012}, the authors propose an inexact decomposition algorithm for solving distributed optimization problems by employing smoothing techniques and an excessive gap condition. In \cite{lin_universal_2015}, the authors proposed an inexact optimization algorithm with an accelerating strategy. The algorithm permits inexact inner-loop solutions. Sufficient conditions on the inexact inner-loop solutions for convergence are shown for different assumptions on the optimization problem.

In \cite{schmidt_convergence_2011}, an inexact proximal-gradient method, as well as its accelerated version, are introduced. The proximal gradient method, also known as the iterative shrinkage-thresholding algorithm (ISTA) \cite{beck_fISTA_2009}, has two main steps: the first one is to compute the gradient of the smooth objective and the second one is to solve the proximal minimization. The conceptual idea of the inexact proximal-gradient method is to allow errors in these two steps, i.e. the error in the calculation of the gradient and the error in the proximal minimization. The results in \cite{schmidt_convergence_2011} show convergence properties of the inexact proximal-gradient method and provide conditions on the errors, under which convergence of the algorithm can be guaranteed. 

Building on the results in \cite{schmidt_convergence_2011}, we propose two new inexact splitting algorithms, the inexact Alternating Minimization Algorithm (inexact AMA) and its accelerated variant, inexact Fast Alternating Minimization Algorithm (inexact FAMA). The inexact FAMA has been studied in \cite{pu_inexact_2014}, and is expanded in this paper. The contributions of this work are the following:

\begin{itemize}
\item \color{black} We propose the inexact AMA and inexact FAMA algorithms, which are inexact variants of the splitting methods, AMA and FAMA in \cite{tseng_applications_1991} and \cite{goldstein_fast_2012}. We show that applying inexact AMA and inexact FAMA to the primal problem is equivalent to applying the inexact proximal-gradient method (inexact PGM) and the inexact accelerated proximal-gradient method (inexact APGM) in \cite{schmidt_convergence_2011} to the dual problem. Based on this fact, we extend the results in \cite{schmidt_convergence_2011}, and show the convergence properties of inexact AMA and inexact FAMA. We derive complexity upper bounds on the number of iterations to achieve a certain accuracy for the algorithms. By exploiting these complexity upper-bounds, we present sufficient conditions on the errors for convergence of the algorithms. 

\item We study the convergence of the algorithms under bounded errors that do not satisfy the sufficient conditions for convergence and show the complexity upper bounds on the number of iterations for this special case.

\item We apply inexact AMA and inexact FAMA for solving distributed optimization problems with local computational errors. We present the complexity upper bounds of the algorithms for this special case, and show sufficient conditions on the local computational errors for convergence. 

\item We study the special case of quadratic local objective functions, relating to a standard form of distributed MPC problems. We show that if the local quadratic functions are positive definite, then the algorithms converge to the optimal solution with a linear rate. We propose to use the proximal gradient method to solve the local problems. By exploiting the sufficient condition on the local computational errors for the convergence together with a  warm-starting strategy, we provide an approach to certify the number of iterations for the proximal gradient method to solve the local problems to the accuracy required for convergence of the distributed algorithm. The proposed on-line certification method only requires on-line local information. 

\item We demonstrate the performance and the theoretical results for inexact algorithms by solving a randomly generated example of a distributed MPC problem with $40$ subsystems. 
\end{itemize}

\section{Preliminaries}\label{sec:Preliminaries}
\subsection{Notation}
 Let $v\in\mathbb{R}^{n_v}$ be a vector. $\|v\|$ denotes the $l_2$ norm of $v$. Let $\mathbb{C}$ be a subset of $\mathbb{R}^{n_v}$. The projection of any point $v\in\mathbb{R}^{n_v}$ onto the set $\mathbb{C}$ is denoted by $\mbox{Proj}_{\mathbb{C}}(v):=\mbox{argmin}_{w\in\mathbb{C}}\;\|w-v\|$. Let $f: \Theta \rightarrow \Omega$ be a function. The conjugate function of $f$ is defined as $f^{\star}(v) = \sup_{w\in \Theta} (v^{T}w - f(w))$. For a conjugate function, it holds that $q\in \partial f(p) \Leftrightarrow p\in \partial f^{\star}(q)$, where $\partial (\cdot)$ denotes the set of sub-gradients of a function at a given point. Let $f$ be a strongly convex function. $\sigma_f$ denotes the convexity modulus $\left< p-q, v-w\right> \geq \sigma_f \|v-w\|^{2}$, where $p\in \partial f(v)$ and $q\in \partial f(w)$, $\forall v, w\in \Theta$. $L(f)$ denotes a Lipschitz constant of the function $f$, i.e. $\|f(v)-f(w)\|\leq L(f)\|v - w \|$, $\forall v, w\in \Theta$. Let $C$ be a matrix. $\rho(C)$ denotes the $l_2$ norm of the matrix $C^{T}C$. The proximity operator is defined as
 
\begin{equation}\label{eq:proximity operator}
\mbox{prox}_{f}(v) = \mbox{argmin}_{w} \quad f(w) + \frac{1}{2}\|w-v\|^{2}\enspace .
\end{equation}

We note the following equivalence:
\begin{equation}\label{eq:property of the conjugate function}
w^{\star}=\mbox{prox}_{f}(v) \Longleftrightarrow v-w^{\star} \in \partial f(w^{\star})
\end{equation}

We refer to \cite{Bertsekas_Convex_2003} and \cite{bauschke_convex_2011} for details on the definitions and properties above. In this paper, $\tilde{\cdot}$ is used to denote an inexact solution of an optimization problem. The proximity operator with an extra subscript $\epsilon$, i.e. $\tilde{x}=\mbox{prox}_{f,\epsilon}(y)$, means that a maximum computation error $\epsilon$ is allowed in the proximal objective function:
\begin{equation}\label{eq:epsilon error in proximity operator}
f(\tilde{w}) + \frac{1}{2}\|\tilde{w}-v\|^{2} \; \leq\; \epsilon + \mbox{min}_{w} \left\lbrace f(w) + \frac{1}{2}\|w-v\|^{2}\right\rbrace
\end{equation}

\subsection{Inexact Proximal-Gradient Method}\label{se:inexact PGM}
In this section, we will introduce the inexact proximal-gradient method (inexact PGM) proposed in \cite{schmidt_convergence_2011}. It addresses optimization problems of the form given in Problem~\ref{pr:problem ISTA} and requires Assumption~\ref{as:ISTA} for convergence, and Assumption~\ref{as:ISTA with strong convexity} for linear convergence. The algorithm is presented in Algorithm~\ref{al:inexact ISTA}. 
\begin{problem}\label{pr:problem ISTA}
\begin{align*}
\min_{w\in \mathbb R^{n_w}} \quad \Phi(w) =\phi(w)+\psi(w) \enspace .
\end{align*}
\end{problem}
\begin{assumption}\label{as:ISTA}
\begin{itemize}
\item $\phi$ is a convex function with Lipschitz continuous gradient with Lipschitz constant $L(\nabla \phi)$
\item $\psi$ is a lower semi-continuous convex function, not necessarily smooth.
\end{itemize}
\end{assumption}

\begin{assumption}\label{as:ISTA with strong convexity}
\begin{itemize}
\item $\phi$ is a strongly convex function with Lipschitz continuous gradient with a convexity modulus $\sigma_{\phi}$.
\item $\psi$ is a lower semi-continuous convex function, not necessarily smooth.
\end{itemize}
\end{assumption}

\begin{scriptsize}
\begin{algorithm}
\caption{Inexact Proximal-Gradient Method}
\begin{algorithmic}
\REQUIRE Require $\tilde{w}^{0} \in \mathbb R^{n_x}$ and $\tau<\frac{1}{L(\nabla \phi)}$
\FOR {$k=1,2,\cdots$}
\STATE 1: $\tilde{w}^{k} =  \mbox{prox}_{\tau \psi, \epsilon^{k}}(\tilde{w}^{k-1} - \tau (\nabla \phi(\tilde{w}^{k-1}) + e^k))$
\ENDFOR
\end{algorithmic}
\label{al:inexact ISTA}
\end{algorithm}\begin{footnotesize}
\end{footnotesize}
\end{scriptsize}

Inexact PGM in Algorithm~\ref{al:inexact ISTA} allows two kinds of errors: $\{e^{k}\}$ represents the error in the gradient calculations of $\phi$, and $\{\epsilon^{k}\}$ represents the error in the computation of the proximal minimization in (\ref{eq:epsilon error in proximity operator}) at every iteration $k$. The following propositions state the convergence property of inexact PGM with different assumptions.
%
%
\begin{proposition}[Proposition 1 in \cite{schmidt_convergence_2011}]\label{pr:convergence rate of inexact ISTA}
Let $\{\tilde{w}_k\}$ be generated by inexact PGM defined in Algorithm~\ref{al:inexact ISTA}. If Assumption \ref{as:ISTA} holds, then for any $k\geq 0$ we have:
\begin{align*}
\Phi\left(\frac{1}{k}\sum^{k}_{p=1} \tilde{w}^p\right)-\Phi(w^{\star})\leq\frac{L(\nabla \phi)}{2k}\left(\|\tilde{w}^{0}-x^{\star}\|+2\Gamma^{k}+\sqrt{2\Lambda^{k}}\right)^{2}
\end{align*}
where $\Phi(\cdot)$ is defined in Problem~\ref{pr:problem ISTA}, 
\begin{equation*}
\Gamma^{k} = \sum^{k}_{p=1} \left(\frac{\|e^{p}\|}{L(\nabla \phi)} +\sqrt{\frac{2\epsilon^{p}}{L(\nabla \phi)}}\right),\;\Lambda^{k} = \sum^{k}_{p=1} \frac{\epsilon^{p}}{L(\nabla \phi)},
\end{equation*}
and $\tilde{w}^{0}$ and $w^{\star}$ denote the initial sequences of Algorithm~\ref{al:inexact ISTA} and the optimal solution of Problem~\ref{pr:problem ISTA}, respectively. 
\end{proposition}
As discussed in \cite{schmidt_convergence_2011}, the complexity upper-bound in Proposition~\ref{pr:convergence rate of inexact ISTA} allows one to derive sufficient conditions on the error sequences $\{e^{k}\}$ and $\{\epsilon^{k}\}$ for the convergence of the algorithm to the optimal solution $w^*$:
\begin{itemize}
\item The series $\{\|e^{k}\|\}$ and $\{\sqrt{\epsilon^{k}}\}$ are finitely summable, i.e., $\sum^{\infty}_{k=1} \|e^{k}\| < \infty$ and $\sum^{\infty}_{k=0} \sqrt{\epsilon^{k}} < \infty$.
\item The sequences $\{\|e^{k}\|\}$ and $\{\sqrt{\epsilon^{k}}\}$ decrease at the rate $O(\frac{1}{k^{1+\kappa}})$ for any $\kappa \geq 0$. 
\end{itemize}

\begin{proposition}[Proposition 3 in \cite{schmidt_convergence_2011}]\label{pr:convergence rate of inexact ISTA with strong convexity}
Let $\{x^k\}$ be generated by inexact PGM defined in Algorithm~\ref{al:inexact ISTA}. If Assumption \ref{as:ISTA with strong convexity} holds, then for any $k\geq 0$ we have:
\begin{align}
\|\tilde{w}^k-w^{\star}\|\leq (1-\gamma)^{k}\cdot(\|w^{0}-w^{\star}\|+\Gamma^{k})\enspace ,
\end{align}
where $\gamma = \frac{\sigma_{\phi}}{L(\nabla \phi)}$ and $w^{0}$ and $w^{\star}$ denote the initial sequence of  Algorithm~\ref{al:inexact ISTA} and the optimal solution of Problem~\ref{pr:problem ISTA}, respectively, and 
\begin{equation*}
\Gamma^{k} = \sum^{k}_{p=1} (1-\gamma)^{-p}\cdot\left(\frac{1}{L(\nabla \phi)} \|e^{p}\| +\sqrt{\frac{2}{L(\nabla \phi)}}\sqrt{\epsilon^{p}}\right)\enspace .
\end{equation*} 
\end{proposition}
From the discussion in \cite{schmidt_convergence_2011}, we can conclude that, if the series $\{\|e^{k}\|\}$ and $\{\sqrt{\epsilon^{k}}\}$ decrease at a linear rate, then $\|x^{k}-x^{\star}\|$ converges to the optimal solution.
%
%
\subsection{Inexact Accelerated Proximal-Gradient Method}\label{se:inexact APGM}
 In this section, we introduce an accelerated variant of inexact PGM, named the inexact accelerated proximal-gradient method (inexact APGM) proposed in \cite{schmidt_convergence_2011}. It addresses the same problem class in Problem~\ref{pr:problem ISTA}  and similarly requires Assumption~\ref{as:ISTA} for convergence. 

%
%
\begin{scriptsize}
\begin{algorithm}
\caption{Inexact Accelerated Proximal-Gradient Method}
\begin{algorithmic}
\REQUIRE Initialize $v^{1} = \tilde{w}^{0} \in \mathbb R^{n_w}$ and $\tau<\frac{1}{L(\nabla \phi)}$
\FOR {$k=1,2,\cdots$} 
\STATE 1: $\tilde{w}^{k} =  \mbox{prox}_{\tau \psi, \epsilon^{k}}(v^{k-1} - \tau (\nabla \phi(v^{k-1}) + e^k))$
\STATE 2: $v^{k} = \tilde{w}^k + \frac{k-1}{k+2}(\tilde{w}^{k}-\tilde{w}^{k-1})$
\ENDFOR
\end{algorithmic}
\label{al:inexact FISTA}
\end{algorithm}\begin{footnotesize}
\end{footnotesize}
\end{scriptsize}
Differing from inexact PGM, inexact APGM involves one extra linear update in Algorithm 2. If Assumption~\ref{as:ISTA} holds, it improves the convergence rate of the complexity upper-bound from $O(\frac{1}{k})$ to $O(\frac{1}{k^{2}})$. The following proposition states the convergence property of inexact APGM.
%
%
\begin{proposition}[Proposition 2 in \cite{schmidt_convergence_2011}]\label{pr:convergence rate of inexact FISTA}
Let $\{\tilde{w}_k\}$ be generated by inexact APGM defined in Algorithm~\ref{al:inexact ISTA}. If Assumption~\ref{as:ISTA} holds, then for any $k\geq 1$ we have:
\begin{align*}
\Phi(\tilde{w}^k)-\Phi(w^{\star})\leq\frac{2L(\nabla \phi)}{(k+1)^{2}}\left(\|\tilde{w}^{0}-x^{\star}\|+2\Gamma^{k}+\sqrt{2\Lambda^{k}}\right)^{2}
\end{align*}
where $\Phi(\cdot)$ is defined in Problem~\ref{pr:problem ISTA}
\begin{equation*}
\Gamma^{k} = \sum^{k}_{p=1} p\left(\frac{\|e^{p}\|}{L(\nabla \phi)} +\sqrt{\frac{2\epsilon^{p}}{L(\nabla \phi)}}\right),\;\Lambda^{k} = \sum^{k}_{p=1} \frac{p^{2}\epsilon^{p}}{L(\nabla \phi)},
\end{equation*}
and $\tilde{w}^{0}$ and $w^{\star}$ denote the starting sequence of  Algorithm~\ref{al:inexact FISTA} and the optimal solution of Problem~\ref{pr:problem ISTA}, respectively.  
\end{proposition}
The complexity upper-bound in Proposition~\ref{pr:convergence rate of inexact FISTA} provides similar sufficient conditions on the error sequences $\{e^{k}\}$ and $\{\epsilon^{k}\}$ for the convergence of Algorithm~\ref{al:inexact FISTA}:
\begin{itemize}
\item The series $\{k\|e^{k}\|\}$ and $\{k\sqrt{\epsilon^{k}}\}$ are finite summable.
\item The sequences $\{\|e^{k}\|\}$ and $\{\sqrt{\epsilon^{k}}\}$ decrease at the rate $O(\frac{1}{k^{2+\kappa}})$ for $\kappa \geq 0$. 
\end{itemize}



\section{Inexact alternating minimization algorithm and its accelerated variant}\label{sec:IAMA and IFAMA}
%

The inexact proximal gradient method, as well as its accelerated version is limited to the case where both objectives are a function of the same variable. However, many optimization problems from engineering fields, e.g. optimal control and machine learning \cite{boyd_distributed_2011}, are not of this problem type. In order to generalize the problem formulation, we employ the alternating minimization algorithm (AMA) and its accelerated variant in \cite{tseng_applications_1991} and \cite{goldstein_fast_2012}, which cover optimization problems of the form of Problem~\ref{pr:AMA and FAMA problem}. In this section, we extend AMA and its accelerated variant to the inexact case and present the theoretical convergence properties.
\begin{problem}\label{pr:AMA and FAMA problem}
\begin{align*}
&\min \quad f(x)+g(z)\\
&\mbox{s.t.} \quad Ax+Bz = c\nonumber
\end{align*}
\end{problem}
\noindent with variables $x\in \mathbb R^{n_x}$ and $z\in \mathbb R^{n_z}$, where $A\in \mathbb R^{n_c\times n_x}$, $B\in \mathbb R^{n_c\times n_z}$ and $c\in \mathbb R^{n_c}$. $f:\mathbb{R}^{n_x} \rightarrow \mathbb{R}$ and $g:\mathbb{R}^{n_z} \rightarrow \mathbb{R}$ are convex functions. The Lagrangian of Problem~\ref{pr:AMA and FAMA problem} is:
\begin{equation}\label{eq:Lagrangian of FAMA problem}
L(x,z,\lambda) = f(x) + g(z) - \lambda^{T}(Ax + Bz -c)\enspace ,
\end{equation}
and the dual function is:
\begin{subequations}
\begin{align}\label{eq:dual function of FAMA problem}
D(\lambda) &= \;\inf_{x,z}\;L(x,z,\lambda) \\
&= -\sup_{x} \left\lbrace \lambda^{T}Ax-f(x)\right\rbrace \nonumber  - \sup_{z} \left\lbrace \lambda^{T}Bz-g(z)\right\rbrace + \lambda^{T}c\\
& = -f^{\star}(A^{T}\lambda) - g^{\star}(B^{T}\lambda) + \lambda^{T}c,
\end{align}
\end{subequations}
where $f^{\star}$ and $g^{\star}$ are the conjugate functions of $f$ and $g$. The dual problem of Problem \ref{pr:AMA and FAMA problem}  is:
\begin{problem}\label{pr:dual problem of AMA and FAMA problem}
\begin{equation*}
\min \quad -D(\lambda) = \underbrace{f^{\star}(A^{T}\lambda)}_{\phi(\lambda)}  + \underbrace{g^{\star}(B^{T}\lambda) - c^{T}\lambda}_{\psi(\lambda)}.
\end{equation*}
\end{problem}

\subsection{Inexact alternating minimization algorithm (inexact AMA)}\label{sec: IAMA}

We propose the inexact alternating minimization algorithm (inexact AMA) presented in Algorithm~\ref{al:Inexact AMA} for solving Problem~\ref{pr:AMA and FAMA problem}. The algorithm allows errors in Step~1 and Step~2, i.e. both minimization problems are solved inexactly with errors $\delta^{k}$ and $\theta^{k}$, respectively. 

\begin{scriptsize}
\begin{algorithm}
\caption{Inexact alternating minimization algorithm (Inexact AMA)}
\begin{algorithmic} 
\REQUIRE Initialize $\lambda^{0}\in \mathbb R^{N_{b}}$, and $\tau < \sigma_{f}/\rho(A)$
\FOR {$k=1,2,\cdots$} 
\STATE 1: $\tilde{x}^{k} = \mbox{argmin}_{x} \;\lbrace f(x) + \langle \lambda^{k-1}, -Ax\rangle\rbrace + \delta^{k} $.
\STATE 2: $\tilde{z}^{k}= \mbox{argmin}_{z} \;\lbrace g(z) + \langle \lambda^{k-1},-Bz\rangle + \frac{\tau}{2}\| c-A\tilde{x}^{k}- Bz\|^{2}\rbrace +\theta^{k}$
\STATE 3: $\lambda^{k} = \lambda^{k-1} + \tau(c-A\tilde{x}^{k}-B\tilde{z}^{k})$
\ENDFOR
\end{algorithmic}
\label{al:Inexact AMA}
\end{algorithm}\begin{footnotesize}
\end{footnotesize}
\end{scriptsize}

We study the theoretical properties of inexact AMA under Assumption~\ref{as:AMA and FAMA}. If Assumption~\ref{as:AMA and FAMA} holds, we show that inexact AMA in Algorithm~\ref{al:Inexact AMA} is equivalent to applying inexact PGM to the dual problem in Problem~\ref{pr:dual problem of AMA and FAMA problem} with the following correspondence: the gradient computation error in Algorithm~\ref{al:inexact FISTA} is equal to $e^{k}=A\delta^{k}$ and the error of solving the proximal minimization is equal to $\epsilon^{k}=\tau ^{2}L(\psi)\|B\theta^{k}\| + \frac{\tau^{2}}{2}\|B\theta^{k}\|^{2}$. With this equivalence, the complexity bound in Proposition~\ref{pr:convergence rate of inexact ISTA} can be extended to the inexact AMA algorithm in Theorem~\ref{th:convergence rate of inexact AMA}. 


%
\begin{assumption}\label{as:AMA and FAMA} We assume that 
\begin{itemize}
\item $f$ is a strongly convex function with convexity modulus $\sigma_{f}$,
\item and $g$ is a convex function, not necessarily smooth.
\end{itemize}
\end{assumption}

\begin{lemma}\label{le:the equivalence between Inexact AMA and inexact ISTA}
If Assumption~\ref{as:AMA and FAMA} is satisfied and inexact AMA and inexact PGM are initialized with the same dual and primal starting sequence, then applying the inexact AMA in Algorithm~\ref{al:Inexact AMA} to Problem~\ref{pr:AMA and FAMA problem} is equivalent to applying inexact PGM in Algorithm~\ref{al:inexact ISTA} to the dual problem defined in Problem~\ref{pr:dual problem of AMA and FAMA problem} with the errors $e^{k} = A\delta^{k}$ and $\epsilon^{k}=\tau ^{2}L(\psi)\|B\theta^{k}\| + \frac{\tau^{2}}{2}\|B\theta^{k}\|^{2}$, where $L(\psi)$ denotes the Lipschitz constant of the function $\psi$.
\end{lemma}

 The proof of Lemma~\ref{le:the equivalence between Inexact AMA and inexact ISTA} is provided in the appendix in Section~\ref{SEC: appendix proof of le:the equivalence between Inexact AMA and inexact ISTA}. This proof is an extension of the proof of Theorem~2 in \cite{goldstein_fast_2012} and the proof in Section 3 in \cite{tseng_applications_1991}. Based on the equivalence shown in Lemma~\ref{le:the equivalence between Inexact AMA and inexact ISTA} , we can now derive an upper-bound on the difference of the dual function value of the sequence $\{\lambda^k\}$ from the optimal dual function value in Theorem~\ref{th:convergence rate of inexact AMA}.
\begin{theorem}\label{th:convergence rate of inexact AMA}
Let $\{\lambda^k\}$ be generated by the inexact AMA in Algorithm~\ref{al:Inexact AMA}. If Assumption~\ref{as:AMA and FAMA} holds, then for any $k\geq 1$
\begin{align}
D(\lambda^{\star}) - D\left(\frac{1}{k}\sum^{k}_{p=1} \lambda^p \right) \leq \frac{L(\nabla \phi)}{2k}\left(\|\lambda^{0}-\lambda^{\star}\|+2\Gamma^{k}+\sqrt{2\Lambda^{k}}\right)^{2} \label{eq:upper bound in the theorem of inexact AMA}
\end{align} 
where $L(\nabla \phi) = \sigma_{f}^{-1}\cdot\rho(A)$,
\begin{align}
&\Gamma^{k} = \sum^{k}_{p=1} \left(\frac{\|A\delta^{p}\|}{L(\nabla \phi)} +\tau\sqrt{\frac{2L(\psi)\|B\theta^{p}\| + \|B\theta^{p}\|^{2}}{L(\nabla \phi)}}\right),\label{eq: Gamma in the bound in the theorem of inexact AMA}\\
& \Lambda^{k} = \sum^{k}_{p=1} \frac{\tau ^{2}(2L(\psi)\|B\theta^{p}\| + \|B\theta^{p}\|^{2})}{2L(\nabla \phi)}\label{eq: Lambda in the bound in the theorem of inexact AMA}
\end{align}
and $\lambda^{0}$ and $\lambda^{\star}$ denote the initial sequences of Algorithm~\ref{al:Inexact AMA} and the optimal solution of Problem~\ref{pr:AMA and FAMA problem}, respectively.
\end{theorem}

\begin{IEEEproof}
Lemma~\ref{le:the equivalence between Inexact AMA and inexact ISTA} shows the equivalence between Algorithm~\ref{al:Inexact AMA} and Algorithm~\ref{al:inexact ISTA} with $e^{k} = A\delta^{k}$ and $\epsilon^{k}=\tau ^{2}L(\psi)\|B\theta^{k}\| + \frac{\tau^{2}}{2}\|B\theta^{k}\|^{2}$. Then we need to show that the dual defined in Problem~\ref{pr:dual problem of AMA and FAMA problem} satisfies Assumption~\ref{as:ISTA}. $\phi(\lambda)$ and $\psi(\lambda)$ are both convex, since the conjugate functions and linear functions as well as their weighted sum are always convex (the conjugate function is the point-wise supremum of a set of affine functions). Furthermore, since $f(x)$ is strongly convex with $\sigma_f$ by Assumption~\ref{as:AMA and FAMA}, then we know $f^{\star}$ has Lipschitz-continuous gradient with Lipschitz constant:
\begin{equation*}\label{eq:lipschitzconstant}
L(\nabla f^{\star}) = \sigma_{f}^{-1}.
\end{equation*}
It follows that the function $\phi$ has Lipschitz-continuous gradient $\nabla \phi$ with a Lipschitz constant:
\begin{equation*}\label{eq:lipschitz constant of phi}
L(\nabla \phi) = \sigma_{f}^{-1}\cdot\rho(A).
\end{equation*}
Hence, the functions $\phi$ and $\psi$ satisfy Assumption~\ref{as:ISTA}. Proposition~\ref{pr:convergence rate of inexact ISTA} completes the proof of the upper-bound in (\ref{eq:upper bound in the theorem of inexact AMA}). 
\end{IEEEproof}

 \color{black}
 
Using the complexity upper-bound in Theorem~\ref{th:convergence rate of inexact AMA}, we derive sufficient conditions on the error sequences for the convergence of inexact AMA in Corollary~\ref{co:sufficient conditions on the errors for IAMA}.

\begin{corollary}\label{co:sufficient conditions on the errors for IAMA}
Let $\{\lambda^k\}$ be generated by the inexact AMA in Algorithm~\ref{al:Inexact AMA}. If Assumption~\ref{as:AMA and FAMA} holds, and the constant $L(\psi) < \infty$, the following sufficient conditions on the error sequences $\{\delta^{k}\}$ and $\{\theta^{k}\}$ guarantee the convergence of Algorithm~\ref{al:Inexact AMA}:
\begin{itemize}
\item The sequences $\{\|\delta^{k}\|\}$ and $\{\|\theta^{k}\|\}$ are finitely summable, i.e., $\sum^{\infty}_{k=1} \|\delta^{k}\| < \infty$ and $\sum^{\infty}_{k=0} \|\theta^{k}\| < \infty$.
\item The sequences $\{\|\delta^{k}\|\}$ and $\{\|\theta^{k}\|\}$ decrease at the rate $O(\frac{1}{k^{1+\kappa}})$ for any $\kappa > 0$. 
\end{itemize}
\end{corollary}

\begin{IEEEproof}
By Assumption~\ref{as:AMA and FAMA}, the dual Problem~\ref{pr:dual problem of AMA and FAMA problem} satisfies Assumption~\ref{al:inexact ISTA} and the complexity upper-bound in Proposition~\ref{pr:convergence rate of inexact ISTA} holds. By extending the sufficient conditions on the error sequences for the convergence of inexact proximal-gradient method discussed after Proposition~\ref{pr:convergence rate of inexact ISTA}, we can derive sufficient conditions on the error sequences for inexact AMA with the errors defined in Lemma~\ref{le:the equivalence between Inexact AMA and inexact ISTA} $e^{k} = A\delta^{k}$ and $\epsilon^{k}=\tau ^{2}L(\psi)\|B\theta^{k}\| + \frac{\tau^{2}}{2}\|B\theta^{k}\|^{2}$. Since $L(\psi) < \infty$, we have that if the error sequences $\{\|\delta^{k}\|\}$ and $\{\|\theta^{k}\|\}$ satisfy the conditions in Corollary~\ref{co:sufficient conditions on the errors for IAMA}, the complexity upper-bound in Theorem~\ref{th:convergence rate of inexact AMA} converges to zero, as the number of iterations $k$ goes to infinity, which further implies that the inexact AMA algorithm converges to the optimal solution.
\end{IEEEproof}

\begin{remark}\label{re:sufficient conditions on the errors for IAMA with psi is an indicator function}
If the function $\psi$ is an indicator function on a convex set, then the constant $L(\psi)$ is equal to infinity, if for any iteration the inexact solution is infeasible with respect to the convex set. However, if we can guarantee that for every iteration $k$ the solutions are feasible with respect to the convex set, then the constant $L(\psi)$ is equal to zero.
\end{remark}

\color{black}

\subsubsection{Linear convergence of inexact AMA for a quadratic cost}\label{sec: IFAMA with linear convergence}

In this section, we study the convergence properties of inexact AMA with a stronger assumption, i.e. the first objective $f$ is a quadratic function and coupling matrix $A$ has full-row rank. We show that with this stronger assumption, the convergence rate of inexact AMA is improved to be linear. The applications satisfying this assumption include least squares problems and distributed MPC problems. 

%
\begin{assumption}\label{as: AMA with quadratic f and A=I}
We assume that
\begin{itemize}
\item $f$ is a quadratic function $f = x^{T}Hx + h^{T}x$ with $H\succ 0$,
\item $A$ has full-row rank.
\end{itemize}
\end{assumption}

\begin{remark}
If Assumption~\ref{as: AMA with quadratic f and A=I} holds, we know that the first objective $\phi(\lambda)$ in the dual problem in Problem~\ref{pr:dual problem of AMA and FAMA problem} is equal to $\phi(\lambda)= \frac{1}{4}(A^{T}\lambda -h)^{T}H^{-1}(A^{T}\lambda -h)$. Then, a Lipschitz constant $L(\nabla \phi)$ is given by the largest eigenvalue of the matrix $\frac{1}{4}AH^{-1}A^{T}$, i.e., $L(\nabla \phi) = \mathbf{\lambda}_{max}(\frac{1}{4}AH^{-1}A^{T})$. In addition, the convexity modulus of $\phi(\lambda)$ is equal to the smallest eigenvalue, i.e., $\sigma_{\phi} = \mathbf{\lambda}_{min}(\frac{1}{4}AH^{-1}A^{T})$. 
\end{remark}

\begin{theorem}\label{th:linear convergence rate of inexact AMA}
Let $\{\lambda^k\}$ be generated by inexact AMA in Algorithm~\ref{al:Inexact AMA}. If Assumption~\ref{as:AMA and FAMA} and ~\ref{as: AMA with quadratic f and A=I} hold, then for any $k\geq 1$
\begin{align}
\|\lambda^k - \lambda^{\star}\|\leq (1-\gamma)^{k}\cdot(\|\lambda^{0}-\lambda^{\star}\|+\Gamma^{k})\enspace , \label{eq:upper bound in the theorem of inexact AMA with linear rate}
\end{align} 
with 
\begin{equation*}
\gamma = \frac{\mathbf{\lambda}_{min}(AH^{-1}A^{T})}{\mathbf{\lambda}_{max}(AH^{-1}A^{T})},
\end{equation*}
%
%
\begin{align*}
\Gamma^{k} = \sum^{k}_{p=1} (1-\gamma)^{-p}\cdot\left(\frac{\|A\delta^{p}\|}{L(\nabla \phi)}  + \tau\sqrt{\frac{L(\psi)\|B\theta^{p}\| + \|B\theta^{p}\|^{2}}{L(\nabla \phi)}}\right).
\end{align*} 
and $\lambda^{0}$ and $\lambda^{\star}$ denote the initial sequences of Algorithm~\ref{al:Inexact AMA} and the optimal solution of Problem~\ref{pr:AMA and FAMA problem}, respectively.
\end{theorem}

\begin{IEEEproof}
By Assumption~\ref{as:AMA and FAMA} and ~\ref{as: AMA with quadratic f and A=I}, the dual problem in Problem~\ref{pr:dual problem of AMA and FAMA problem} satisfies Assumption~\ref{as:ISTA with strong convexity} and the complexity upper-bound in Proposition~\ref{pr:convergence rate of inexact ISTA with strong convexity} holds for the dual problem. The proof of Theorem~\ref{th:linear convergence rate of inexact AMA} follows directly from this fact.
\end{IEEEproof}

\color{black}	
By using the complexity upper-bounds in Theorem~\ref{th:linear convergence rate of inexact AMA}, we derive sufficient conditions on the error sequences, which guarantee the convergence of the inexact AMA algorithm.

\begin{corollary}\label{co:sufficient conditions on the errors for IAMA with linear convergence rate}
Let $\{\lambda^k\}$ be generated by the inexact AMA in Algorithm~\ref{al:Inexact AMA}. If Assumption~\ref{as:AMA and FAMA} and ~\ref{as: AMA with quadratic f and A=I} hold, and the constant $L(\psi) < \infty$, the following sufficient conditions on the error sequences $\{\delta^{k}\}$ and $\{\theta^{k}\}$ guarantee the convergence of Algorithm~\ref{al:Inexact AMA}:
\begin{itemize}
\item The sequences $\{\|\delta^{k}\|\}$ and $\{\|\theta^{k}\|\}$ are finitely summable, i.e., $\sum^{\infty}_{k=1} \|\delta^{k}\| < \infty$ and $\sum^{\infty}_{k=0} \|\theta^{k}\| < \infty$.
\item The sequences $\{\|\delta^{k}\|\}$ and $\{\theta^{k}\}$ decrease at  $O(\frac{1}{k^{1+\kappa}})$ for any $\kappa \in \mathbb{Z}_{+}$. For this case the complexity upper-bound in (\ref{eq:upper bound in the theorem of inexact AMA with linear rate}) reduces to the same rate as the error sequences.
\item The sequences $\{\|\delta^{k}\|\}$ and $\{\|\theta^{k}\|\}$ decrease at a linear rate. 
\end{itemize}
\end{corollary}

\begin{IEEEproof}
By using the complexity upper-bound in Theorem~\ref{th:linear convergence rate of inexact AMA}, we can derive sufficient conditions on the error sequences for inexact AMA. Since $L(\psi) < \infty$, we have that if the error sequences $\{\|\delta^{k}\|\}$ and $\{\|\theta^{k}\|\}$ satisfy the first and third conditions in Corollary~\ref{co:sufficient conditions on the errors for IAMA with linear convergence rate}, the complexity upper-bound in Theorem~\ref{th:linear convergence rate of inexact AMA} converges to zero, as the number of iterations $k$ goes to infinity, which further implies that the inexact AMA algorithm converges to the optimal solution. For the second sufficient condition in Corollary~\ref{co:sufficient conditions on the errors for IAMA with linear convergence rate}, we provide Lemma~\ref{le: convergence for series for sufficient condition for IAMA with linear rate} to prove that the second sufficient condition guarantees the convergence of the algorithm.
\end{IEEEproof}
\color{black}	


\color{black}	
\begin{lemma}\label{le: convergence for series for sufficient condition for IAMA with linear rate}
Let $\alpha$ be a positive number $0 < \alpha < 1$. The following series $S^{k}$ converges to zero, as the index $k$ goes to infinity
\begin{equation*}
\lim_{k \rightarrow \infty} S^k := \lim_{k \rightarrow \infty} \alpha^{k}\cdot \sum^{k}_{p=1} \frac{\alpha^{-p}}{p}=0\enspace .
\end{equation*}
Furthermore, the series $S^k$ converges at the rate $O(\frac{1}{k})$.
\end{lemma}

 The proof of Lemma~\ref{le: convergence for series for sufficient condition for IAMA with linear rate} is provided in the appendix in Section~\ref{SEC: appendix proof of le: convergence for series for sufficient condition for IAMA with linear rate}. Lemma~\ref{le: convergence for series for sufficient condition for IAMA with linear rate} provides that if the sequences $\{\|\delta^{k}\|\}$ and $\{\|\theta^{k}\|\}$ decrease at $O(\frac{1}{k})$, the complexity upper-bound in (\ref{eq:upper bound in the theorem of inexact AMA with linear rate}) converges at the rate $O(\frac{1}{k})$. Note that this result can be extended to the case that $\{\|\delta^{k}\|\}$ and $\{\|\theta^{k}\|\}$ decrease at  $O(\frac{1}{k^{1+\kappa}})$ for any $\kappa \in \mathbb{Z}_{+}$, by following a similar proof as for Lemma~\ref{le: convergence for series for sufficient condition for IAMA with linear rate}.
 \color{black}

\subsection{Inexact fast alternating minimization algorithm (inexact FAMA)}\label{sec: IFAMA}

 In this section, we present an accelerated variant of inexact AMA, named the inexact fast alternating minimization Algorithm (inexact FAMA), which is presented in Algorithm~\ref{al:Inexact FAMA}. It addresses the same problem class as Problem~\ref{pr:AMA and FAMA problem} and requires the same assumption as Assumption~\ref{as:AMA and FAMA} for convergence. Similar to inexact AMA, inexact FAMA allows computation errors in the two minimization steps in the algorithm. Differing from inexact AMA, inexact FAMA involves one extra linear update in Step~4 in Algorithm~\ref{al:Inexact FAMA}, which improves the optimal convergence rate of the complexity upper-bound of the algorithm from $O(\frac{1}{k})$ to $O(\frac{1}{k^{2}})$. This is similar to the relationship between the inexact PGM and inexact APGM. By extending the result in Lemma~\ref{le:the equivalence between Inexact AMA and inexact ISTA}, we show that inexact FAMA is equivalent to applying inexact APGM to the dual problem. With this equivalence, we further show a complexity upper bound for inexact FAMA by using the result in Proposition~\ref{pr:convergence rate of inexact FISTA} for inexact APGM. \color{black}The main results for inexact FAMA have been presented in \cite{pu_inexact_2014}, and are restated in this section.\color{black}
\begin{scriptsize}
\begin{algorithm}
\caption{Inexact Fast alternating minimization algorithm (Inexact FAMA)}
\begin{algorithmic} 
\REQUIRE Initialize $ \hat{\lambda}^{0} = \lambda^{0}\in \mathbb R^{N_{b}}$, and $\tau < \sigma_{f}/\rho(A)$
\FOR {$k=1,2,\cdots$}
\STATE 1: $\tilde{x}^{k} = \mbox{argmin}_{x} \;\lbrace f(x) + \langle \hat{\lambda}^{k-1}, -Ax\rangle\rbrace + \delta^{k} $.
\STATE 2:$\tilde{z}^{k}= \mbox{argmin}_{z} \;\lbrace g(z) + \langle \hat{\lambda}^{k-1},-Bz\rangle+ \frac{\tau}{2}\| c-A\tilde{x}^{k}$
$ - Bz\|^{2}\rbrace +\theta^{k}$
\STATE 3: $\lambda^{k} = \hat{\lambda}^{k-1} + \tau(c-A\tilde{x}^{k}-B\tilde{z}^{k})$
\STATE 4: $\hat{\lambda}^{k} = \lambda^{k} + \frac{k-1}{k+2}(\lambda^{k}-\lambda^{k-1})$
\ENDFOR
\end{algorithmic}
\label{al:Inexact FAMA}
\end{algorithm}\begin{footnotesize}
\end{footnotesize}
\end{scriptsize}
\begin{lemma}\label{le:the equivalence between Inexact FAMA and inexact FISTA}
If Assumption~\ref{as:AMA and FAMA} is satisfied and inexact FAMA and inexact APGM are initialized with the same dual and primal starting sequence, respectively, then applying the inexact FAMA in Algorithm~\ref{al:Inexact FAMA} to Problem~\ref{pr:AMA and FAMA problem} is equivalent to applying inexact APGM in Algorithm~\ref{al:inexact FISTA} to the dual problem defined in Problem~\ref{pr:dual problem of AMA and FAMA problem} with the errors $e^{k} = A\delta^{k}$ and $\epsilon^{k}=\tau ^{2}L(\psi)\|B\theta^{k}\| + \frac{\tau^{2}}{2}\|B\theta^{k}\|^{2}$, where $L(\psi)$ denotes the Lipschitz constant of the function $\psi$.
\end{lemma}

\begin{IEEEproof}
The proof follows the same flow of the proof of Lemma~\ref{le:the equivalence between Inexact AMA and inexact ISTA} by replacing $\lambda^{k-1}$ by $\hat{\lambda}^{k-1}$ computed in Step~4 in Algorithm~\ref{al:Inexact FAMA} and showing the following equality
\begin{equation}\label{eq: equivalence of inexact FAMA and inexact FISTA on dual}
\lambda^{k} = \mbox{prox}_{\tau \psi, \epsilon^{k}}(\hat{\lambda}^{k-1} - \tau (\nabla \phi(\hat{\lambda}^{k-1}) + e^{k}))
\end{equation}
%
\end{IEEEproof}

Based on the equivalence shown in Lemma~\ref{le:the equivalence between Inexact FAMA and inexact FISTA} , we can now derive an upper-bound on the difference of the dual function value of the sequence $\{\lambda^k\}$ for inexact FAMA in Theorem~\ref{th:convergence rate of inexact FAMA}.
\begin{theorem}[Theorem~III.5 in \cite{pu_inexact_2014}]\label{th:convergence rate of inexact FAMA}
Let $\{\lambda^k\}$ be generated by the inexact FAMA in Algorithm~\ref{al:Inexact FAMA}. If Assumption~\ref{as:AMA and FAMA} holds, then for any $k\geq 1$
\begin{equation}\label{eq:upper bound in the theorem of inexact FAMA}
D(\lambda^{\star}) - D(\lambda^k) \leq \frac{2L(\nabla \phi)}{(k+1)^{2}}\left(\|\lambda^{0}-\lambda^{\star}\|+2\Gamma^{k}+\sqrt{2\Lambda^{k}}\right)^{2}
\end{equation} 
where 
\begin{align}
&\Gamma^{k} = \sum^{k}_{p=1} p\left(\frac{\|A\delta^{p}\|}{L(\nabla \phi)} +\tau\sqrt{\frac{2L(\psi)\|B\theta^{p}\| + \|B\theta^{p}\|^{2}}{L(\nabla \phi)}}\right),\label{eq: Gamma in the bound in the theorem of inexact FAMA}\\
& \Lambda^{k} = \sum^{k}_{p=1} \frac{p^{2}\tau ^{2}(2L(\psi)\|B\theta^{p}\| + \|B\theta^{p}\|^{2})}{2L(\nabla \phi)}\label{eq: Lambda in the bound in the theorem of inexact FAMA}
\end{align}
and $L(\nabla \phi) = \sigma_{f}^{-1}\cdot\rho(A)$.
\end{theorem}
\begin{IEEEproof}
The proof is a similar to the proof of Theorem~\ref{th:convergence rate of inexact AMA}. Lemma~\ref{le:the equivalence between Inexact FAMA and inexact FISTA} shows the equivalence between Algorithm~\ref{al:Inexact FAMA} and Algorithm~\ref{al:inexact FISTA}. Proposition~\ref{pr:convergence rate of inexact FISTA} completes the proof of the upper-bound in inequality (\ref{eq:upper bound in the theorem of inexact FAMA}). 
\end{IEEEproof}

With the results in Theorem~\ref{th:convergence rate of inexact FAMA}, the sufficient conditions on the errors for the convergence of inexact APGM presented in Section~\ref{se:inexact APGM} can be extended to inexact FAMA with the errors defined in Lemma~\ref{le:the equivalence between Inexact FAMA and inexact FISTA}.

\color{black}
\begin{corollary}\label{co:sufficient conditions on the errors for IFAMA}
Let $\{\lambda^k\}$ be generated by the inexact AMA in Algorithm~\ref{al:Inexact FAMA}. If Assumption~\ref{as:AMA and FAMA} holds, and the constant $L(\psi) < \infty$, the following sufficient conditions on the error sequences $\{\delta^{k}\}$ and $\{\theta^{k}\}$ guarantee the convergence of Algorithm~\ref{al:Inexact AMA}:
\begin{itemize}
\item The series $\{\|\delta^{k}\|\}$ and $\{\|\theta^{k}\|\}$ are finitely summable, i.e., $\sum^{\infty}_{k=1} \|\delta^{k}\| < \infty$ and $\sum^{\infty}_{k=0} \|\theta^{k}\| < \infty$.
\item The sequences $\{\|\delta^{k}\|\}$ and $\{\|\theta^{k}\|\}$ decrease at the rate $O(\frac{1}{k^{2+\kappa}})$ for any $\kappa > 0$. 
\end{itemize}
\end{corollary}

\begin{IEEEproof}
By Assumption~\ref{as:AMA and FAMA}, the dual Problem~\ref{pr:dual problem of AMA and FAMA problem} satisfies Assumption~\ref{al:inexact ISTA} and the complexity upper-bound in Proposition~\ref{pr:convergence rate of inexact FISTA} holds. By extending the sufficient conditions on the error sequences discussed after Proposition~\ref{pr:convergence rate of inexact FISTA}, we obtain sufficient conditions on the error sequences for inexact FAMA with the errors defined in Lemma~\ref{le:the equivalence between Inexact FAMA and inexact FISTA}. Since $L(\psi) < \infty$, we have that if the error sequences $\{\|\delta^{k}\|\}$ and $\{\|\theta^{k}\|\}$ satisfy the conditions in Corollary~\ref{co:sufficient conditions on the errors for IFAMA}, the complexity upper-bound in Theorem~\ref{th:convergence rate of inexact FAMA} converges to zero, as the number of iterations $k$ goes to infinity, which further implies that the inexact FAMA algorithm converges to the optimal solution.
\end{IEEEproof}
\color{black}




\subsection{Discussion: inexact AMA and inexact FAMA with bounded errors}

\color{black}
 In this section, we study the special case that the error sequences $\delta^{k}$ and $\theta^{k}$ are bounded by constants. This special case is of particular interest, as it appears in many engineering problems in practice, e.g. quantized distributed computation and distributed optimization with constant local computation errors. Previous work includes \cite{necoara_rate_2013}, where the authors studied the complexity upper-bounds for a distributed optimization algorithm with bounded noise on the solutions of local problems. In this section, we will study errors satisfying Assumption~\ref{as: the errors are bounded for IFAMA} and derive the corresponding complexity upper-bounds for inexact AMA, as well as for inexact FAMA, with different assumptions. We show that if the problem satisfies the stronger assumption in Assumption~\ref{as: AMA with quadratic f and A=I}, i.e. the cost function $f$ is a quadratic function, then the complexity bounds for the inexact algorithms with bounded errors converge to a finite positive value, as $k$ increases. It is important to point out that if only the conditions in Assumption~\ref{as:AMA and FAMA} are satisfied, convergence of the complexity upper bound to a small constant cannot be shown. We present the complexity upper-bound of inexact FAMA in details for this case, and the result can be easily extended to inexact AMA.
 
 %
\begin{assumption}\label{as: the errors are bounded for IFAMA}
We assume that the error sequences $\delta^{k}$ and $\theta^{k}$ are bounded by $\|\delta^{k}\|\leq \bar{\delta}$ and $\|\theta^{k}\|\leq \bar{\theta}$ for all $k\geq 0$, where $\bar{\delta}$ and $\bar{\theta}$ are positive constants. 
\end{assumption}
\begin{corollary}\label{th:convergence rate of inexact AMA with bounded errors}
Let $\{\lambda^k\}$ be generated by the inexact AMA in Algorithm~\ref{al:Inexact AMA}. If Assumption~\ref{as:AMA and FAMA}, \ref{as: AMA with quadratic f and A=I} and \ref{as: the errors are bounded for IFAMA} hold, then for any $k\geq 1$
\begin{align}
\|\lambda^k - \lambda^{\star}\|\leq (1-\gamma)^{k}\cdot\|\lambda^{0}-\lambda^{\star}\| +\Delta\enspace , \label{eq:upper bound in the theorem of inexact AMA with linear rate with bounded errors}
\end{align} 
where $\Delta = \frac{1}{\gamma}\left(\frac{\|A\bar{\delta}\|}{L(\nabla \phi)}  + \tau\sqrt{\frac{L(\psi)\|B\bar{\theta}\| + \|B\bar{\theta}\|^{2}}{L(\nabla \phi)}}\right)$, $\gamma = \frac{\mathbf{\lambda}_{min}(AH^{-1}A^{T})}{\mathbf{\lambda}_{max}(AH^{-1}A^{T})}$ and $\lambda^{0}$ and $\lambda^{\star}$ denote the initial sequences of Algorithm~\ref{al:Inexact AMA} and the optimal solution of Problem~\ref{pr:AMA and FAMA problem}, respectively.
\end{corollary}
\begin{IEEEproof}
Since Assumption~\ref{as:AMA and FAMA} and~\ref{as: AMA with quadratic f and A=I} are satisfied, then the results in Theorem~\ref{th:linear convergence rate of inexact AMA} hold. By Assumption~\ref{as: the errors are bounded for IFAMA}, we know that the error sequences satisfy $\|\delta^{k}\|\leq \bar{\delta}$ and $\|\theta^{k}\|\leq \bar{\theta}$ for all $k\geq 0$. Then the error function $\Gamma^{k}$ in Theorem~\ref{th:linear convergence rate of inexact AMA} is upper-bounded by
\begin{align*}
\Gamma^{k} \leq\sum^{k}_{p=1} (1-\gamma)^{-p}\cdot\left(\frac{\|A\bar{\delta}\|}{L(\nabla \phi)}  + \tau\sqrt{\frac{L(\psi)\|B\bar{\theta}\| + \|B\bar{\theta}\|^{2}}{L(\nabla \phi)}}\right)\enspace .
\end{align*} 
Due to the fact that $0< \gamma <1$ and the property of geometric series, we get
\begin{align*}
(1-\gamma)^{k}\cdot \Gamma^{k} &\leq \sum^{k}_{p=1} (1-\gamma)^{k-p}\cdot\left(\frac{\|A\bar{\delta}\|}{L(\nabla \phi)}  + \tau\sqrt{\frac{L(\psi)\|B\bar{\theta}\| + \|B\bar{\theta}\|^{2}}{L(\nabla \phi)}}\right)\\
&\leq \frac{1-(1-\gamma)^{k}}{\gamma}\cdot\left(\frac{\|A\bar{\delta}\|}{L(\nabla \phi)}  + \tau\sqrt{\frac{L(\psi)\|B\bar{\theta}\| + \|B\bar{\theta}\|^{2}}{L(\nabla \phi)}}\right)\\
&\leq \frac{1}{\gamma}\cdot\left(\frac{\|A\bar{\delta}\|}{L(\nabla \phi)}  + \tau\sqrt{\frac{L(\psi)\|B\bar{\theta}\| + \|B\bar{\theta}\|^{2}}{L(\nabla \phi)}}\right)\enspace .
\end{align*} 
Then the upper-bound in Theorem~\ref{th:linear convergence rate of inexact AMA} implies the upper-bound in~(\ref{eq:upper bound in the theorem of inexact AMA with linear rate with bounded errors}).
\end{IEEEproof}
\begin{remark}
The inexact AMA algorithm with bounded errors satisfying Assumption~\ref{as:AMA and FAMA} and~\ref{as: AMA with quadratic f and A=I} has a constant term $\Delta$ in the complexity upper-bound in (\ref{eq:upper bound in the theorem of inexact AMA with linear rate with bounded errors}). Hence, the complexity bound in (\ref{eq:upper bound in the theorem of inexact AMA with linear rate with bounded errors}) converges to a neighbourhood of the origin with the size of $\Delta$, as $k$ goes to infinity.
\end{remark}
\begin{remark}\label{re:convergence rate of inexact FAMA with bounded errors}
For the inexact FAMA in Algorithm~\ref{al:Inexact FAMA}, if Assumption~\ref{as:AMA and FAMA} and Assumption~\ref{as: the errors are bounded for IFAMA} hold, i.e. the cost is not necessarily quadratic, we can also derive the following complexity upper bound
\begin{equation}\label{eq: complexity bound for IFAMA with bounded errors}
D(\lambda^{\star}) - D(\lambda^k) \leq \left(\frac{2L(\nabla \phi)\|\lambda^{0}-\lambda^{\star}\|}{(k+1)}+k\cdot \Delta\right)^{2}
\end{equation} 
with 
\begin{align*}
\Delta = \frac{\|A\|\cdot\bar{\delta}}{L(\nabla \phi)} +\frac{3\tau}{2}\cdot\sqrt{\frac{(2L(\psi)\|B\|\cdot\bar{\theta} + \|B\|\cdot\bar{\theta}^{2})}{L(\nabla \phi)}}
\end{align*}
and $L(\nabla \phi) = \sigma_{f}^{-1}\cdot\rho(A)$. The proof follows the same flow of the proof for Corollary~\ref{th:convergence rate of inexact AMA with bounded errors} by replacing Theorem~\ref{th:linear convergence rate of inexact AMA} with Theorem~\ref{th:convergence rate of inexact FAMA}. Compared to the FAMA algorithm without errors, we see that the inexact FAMA with bounded errors has one extra term $k\cdot \Delta$ in the complexity upper-bound in (\ref{eq: complexity bound for IFAMA with bounded errors}). Unfortunately, the term $k\cdot \Delta$ increases as $k$ increases. Hence, the complexity bound for the inexact FAMA with bounded errors does not converge, as $k$ goes to infinity.
\end{remark}
\color{black}

\section{Inexact AMA for distributed optimization with an application to distributed MPC}\label{sec: Inexact FAMA for Distributed Optimization with a Distributed  MPC Application}

\subsection{Distributed optimization problem}

In this section, we consider a distributed optimization problem on a network of $M$ sub-systems (nodes). The sub-systems communicate according to a fixed undirected graph $G =(\mathcal{V},\mathcal{E})$. The vertex set $\mathcal{V} = \{1,2,\cdots,M\}$ represents the sub-systems and the edge set $\mathcal{E}\subseteq \mathcal{V}\times \mathcal{V}$ specifies pairs of sub-systems that can communicate. If $(i,j)\in \mathcal{E}$, we say that sub-systems $i$ and $j$ are neighbours, and we denote by $\mathcal{N}_i = \{j| (i,j)\in \mathcal{E}\}$ the set of the neighbours of sub-system $i$. Note that $\mathcal{N}_i$ includes $i$. The cardinality of $\mathcal{N}_i$ is denoted by $|\mathcal{N}_i|$. The global optimization variable is denoted by $z$. The local variable of sub-system $i$, namely the $i$th element of $z$ and $z=[z^{T}_1,\cdots,z^{T}_M]^{T}$, is denoted by $[z]_i$. The concatenation of the variable of sub-system $i$ and the variables of its neighbours is denoted by $z_{i}$. With the selection matrices $E_i$ and $F_{ji}$, the variables have the following relationship: $z_{i} = E_iz$ and $[z]_i = F_{ji}z_{j}$, $j\in\mathcal{N}_i $, which implies the relation between the local variable $[z]_i$ and the global variable $z$, i.e. $[z]_i = F_{ji}E_jz$, $j\in\mathcal{N}_i $. We consider the following distributed optimization problem:
\begin{problem}\label{pr:split form for distributed MPC}
\begin{align*}
\min_{z,v}  & \quad \sum^{M}_{i=1} f_{i}(z_i)\\
s.t. &\quad  z_i\in \mathbb{C}_{i}, \quad z_i = E_iv, \quad i=1,2,\cdots ,M.
\end{align*}
\end{problem}
\noindent where $f_i$ is the local cost function for sub-system $i$, and the constraint $\mathbb{C}_i$ represents a convex local constraint on the concatenation of the variable of sub-system $i$ and the variables of its neighbours $z_{i}$.
\begin{assumption}\label{as: distributed problems for inexact AMA}
Each local cost function $f_i$ in Problem~\ref{pr:split form for distributed MPC} is a strongly convex function with a convexity modulus $\sigma_{f_i}$ and has a Lipschitz continuous gradient with Lipschitz constant $L(\nabla f_i)$. The set $\mathbb{C}_{i}$ is a convex set, for all $i=1,\cdots,M$.
\end{assumption}

\begin{remark}\label{rem:IAMA for DMPC constrains on the first objective}
Recall the problem formulation of inexact AMA and FAMA defined in Problem~\ref{pr:AMA and FAMA problem}. The two objectives are defined as $f(z)=\sum^{M}_{i=1} f_{i}(z_i)$ subject to $z_i\in \mathbb{C}_{i}$ for all $i = 1,\cdots,M$ and $g=0$. The matrices are $A=I$, $B=-[E^{T}_1,E^{T}_2,\cdots ,E^{T}_M]^{T}$ and $c=0$. The first objective $f(\mathbf{z})$ consists of a strongly convex function on $z$ and convex constraints. The convex constraints can be considered as indicator functions, which are convex functions. Due to the fact that the sum of a strongly convex and a convex function is strongly convex, the objective $f(z)$ is strongly convex with the modulus $\sigma_{f}$ and Problem~\ref{pr:split form for distributed MPC} satisfies Assumption~\ref{as:AMA and FAMA}.
\end{remark}

\subsection{Application: distributed model predictive control}\label{se: DMPC}

In this section, we consider a distributed linear MPC problem with $M$ sub-systems, and show that it can be written in the form of Problem~\ref{pr:split form for distributed MPC}. The dynamics of the $i$th agent are given by the discrete-time linear dynamics:
%
%
\begin{equation}
x_{i}(t+1) = \sum_{j\in \mathcal{N}_{i}} A_{ij}x_{j}(t) + B_{ij}u_{j}(t) \quad i=1,2,\cdots ,M.
\end{equation}
where $A_{ij}$ and $B_{ij}$ are the dynamical matrices. The states and inputs of agent $i$ are subject to local convex constraints:
\begin{equation}
x_i(t) \in \mathbb{X}_i\quad u_i(t) \in \mathbb{U}_i\quad i=1,2,\cdots ,M.
\end{equation}
The distributed MPC problem, as e.g. considered in \cite{conte_distributed_2012}, is given in Problem~\ref{pr:distributed MPC}.
\begin{problem}\label{pr:distributed MPC}
\begin{align*}
\min_{x,u}  & \quad \sum^{M}_{i=1} \sum^{N-1}_{t=0} l_i(x_{i}(t),u_{i}(t)) + \sum^{M}_{i=1} l^{f}_i(x_{i}(N)) \\
s.t. \quad &x_{i}(t+1) = \sum_{j\in \mathcal{N}_{i}} A_{ij}x_{j}(t) + B_{ij}u_{j}(t)\\
& x_i(t) \in \mathbb{X}_i\quad u_i(t) \in \mathbb{U}_i\\
& x_i(N) \in \mathbb{X}^{f}_i ,\quad x_i(0) = \bar{x}_{i}, \quad i=1,2,\cdots ,M.
\end{align*}
\end{problem}
\noindent where $l_i(\cdot,\cdot)$ and $l^{f}_i(\cdot)$ are strictly convex stage cost functions and $N$ is the horizon for the MPC problem. The state and input sequences along the horizon of agent $i$ are denoted by $x_{i}=[x^{T}_{i}(0),x^{T}_{i}(1),\cdots,x^{T}_{i}(N)]^{T}$ and $u_{i}=[u^{T}_{i}(0),u^{T}_{i}(1),\cdots,u^{T}_{i}(N-1)]^{T}$. We denote the concatenations of the state and input sequences of agent $i$ and its neighbours by $x_{\mathcal{N}_i}$ and $u_{\mathcal{N}_i}$. The corresponding constraints are $x_{\mathcal{N}_i} \in \mathbb{X}_{\mathcal{N}_i}$ and $u_{\mathcal{N}_i} \in \mathbb{U}_{\mathcal{N}_i}$. We define $v=[x^{T}_{1},x^{T}_{2},\cdots,x^{T}_{M},u^{T}_{1},u^{T}_{2},\cdots,u^{T}_{M}]^{T}$ to be the global variable and $z_{i} = [x_{\mathcal{N}_i},u_{\mathcal{N}_i}]$ to be the local variables. $\mathbb{Z}_{\mathcal{N}_i}=\mathbb{X}_{\mathcal{N}_i}\times\mathbb{U}_{\mathcal{N}_i}$ denotes the local constraints on $z_i$ and $H_iz_i=h_i$ denotes the dynamical constraint of sub-system $i$. Then considering the distributed problem in Problem~\ref{pr:split form for distributed MPC}, we see that the local cost function $f_i$ for agent $i$ contains all the stage cost functions of the state and input sequences of agent $i$ and its neighbours. The constraint $\mathbb{C}_i$ includes the constraint $\mathbb{Z}_{\mathcal{N}_i}$ and the dynamical constraint $H_iz_i=h_i$. $E_i$ are the matrices selecting the local variables from the global variable. The $i$th component of $v$ is equal to $[v]_i = [x_{i},u_{i}]$. 

\begin{remark}\label{re: discussion about quadratic cost functions of DMPC for IAMA}
If the stage cost functions $l_i(\cdot,\cdot)$ and $l^{f}_i(\cdot)$ are strictly convex functions, and the state and input constraints $\mathbb{X}_{i}$ and $\mathbb{U}_{i}$ are convex sets, then the conditions in Assumption~\ref{as: distributed problems for inexact AMA} are all satisfied. Furthermore, if the state cost functions $l_i(\cdot,\cdot)$ and $l^{f}_i(\cdot)$ are set to be positive definite quadratic functions, then the distributed optimization problem originating from the distributed MPC problem further satisfies Assumption~\ref{as: AMA with quadratic f and A=I}.
\end{remark}

\color{black}
\begin{remark}\label{re: discussion condensed form of DMPC (controllable) for IAMA}
For the case that the distributed MPC problem has only input constraints and the state coupling matrices in the linear dynamics are $A_{ij}=0$ for any $i\neq j$, we can eliminate all state variables in the distributed MPC problem and only have the input variables as the optimization variables. For this case, if the stage cost functions $l_i(\cdot,\cdot)$ and $l^{f}_i(\cdot)$ are strictly convex functions with respect to the input variables and the local linear dynamical system $x_{i}(t+1) = A_{ii}x_{i}(t) + \sum_{j\in \mathcal{N}_{i}} B_{ij}u_{j}(t)$ is controllable, then the resulting distributed optimization problem satisfies Assumption~\ref{as: AMA with quadratic f and A=I}. The details of this formulation can be found in \cite{pu_quantization_par_2015}. 
\end{remark}
\color{black}

\subsection{Inexact AMA and Inexact FAMA for distributed optimization}\label{se:IFAMA for DMPC}
In this section, we apply inexact AMA and inexact FAMA to the distributed optimization problem in Problem~\ref{pr:split form for distributed MPC}, originating from the distributed MPC problem in Problem~\ref{pr:distributed MPC}. The concept is to split the distributed optimization into small and local problems according to the physical couplings of the sub-systems. Algorithm~\ref{al:IAMA for DMPC} and Algorithm~\ref{al:IFAMA for DMPC} represent the algorithms. Note that Step~2 in inexact AMA and inexact FAMA, i.e., Algorithm~\ref{al:Inexact AMA} and Algorithm~\ref{al:Inexact FAMA}, are simplified to be a consensus step in Step~3 in Algorithm~\ref{al:IAMA for DMPC} and Algorithm~\ref{al:IFAMA for DMPC}, which requires only local communication. In the algorithms, $\delta^{k}_i$ represents the computational error of the local problems. 

\begin{scriptsize}
\begin{algorithm}
\caption{Inexact Alternating Minimization Algorithm for Distributed Optimization}
\begin{algorithmic} 
\REQUIRE Initialize $\lambda^{0}_i =0 \in \mathbb{R}^{z_i}$, and $\tau < \min_{1\leq i \leq M} \{\sigma_{f_i}\}$  
\FOR {$k=1,2,\cdots$}
\STATE 1: $\tilde{z}^{k}_i = \mbox{argmin}_{z_i\in \mathbb{C}_{i}}\lbrace  f_i(z_i) + \langle \lambda^{k-1}_{i}, -z_i\rangle\rbrace + \delta^{k}_i$
\STATE 2: Send $\tilde{z}^{k}_i$ to all the neighbours of agent $i$.
\STATE 3: $[\tilde{v}^{k}]_{i} = \frac{1}{|\mathcal{N}_i|} \sum^{M}_{j\in\mathcal{N}_i} [\tilde{z}^{k}_j]_i$.
\STATE 4: Send $[\tilde{v}^{k}]_{i}$ to all the neighbours of agent $i$.
\STATE 5: $\lambda^{k}_i = \lambda^{k-1}_i + \tau(E_i\tilde{v}^{k}-\tilde{z}^{k}_i)$
\ENDFOR
\end{algorithmic}
\label{al:IAMA for DMPC}
\end{algorithm}\begin{footnotesize}
\end{footnotesize}
\end{scriptsize}


\begin{scriptsize}
\begin{algorithm}
\caption{Inexact fast alternating minimization algorithm for Distributed Optimization}
\begin{algorithmic} 
\REQUIRE Initialize $\lambda^{0}_i = \hat{\lambda}^{0}_{i} \in \mathbb{R}^{z_i}$, and $\tau < \min_{1\leq i \leq M} \{\sigma_{f_i}\}$
\FOR {$k=1,2,\cdots$}
\STATE 1: $\tilde{z}^{k}_i = \mbox{argmin}_{z_i\in \mathbb{C}_{i}}\lbrace  f_i(z_i) + \langle \hat{\lambda}^{k-1}_{i}, -z_i\rangle\rbrace + \delta^{k}_i$
\STATE 2: Send $\tilde{z}^{k}_i$ to all the neighbours of agent $i$.
\STATE 3: $[\tilde{v}^{k}]_{i} = \frac{1}{|\mathcal{N}_i|} \sum^{M}_{j\in\mathcal{N}_i} [\tilde{z}^{k}_j]_i$.
\STATE 4: Send $[\tilde{v}^{k}]_{i}$ to all the neighbours of agent $i$.
\STATE 5: $\lambda^{k}_i = \hat{\lambda}^{k-1}_i + \tau(E_i\tilde{v}^{k}-\tilde{z}^{k}_i)$
\STATE 6: $\hat{\lambda}^{k}_i = \lambda^k_i + \frac{k-1}{k+2}(\lambda^k_i-\lambda^{k-1}_i)$
\ENDFOR
\end{algorithmic}
\label{al:IFAMA for DMPC}
\end{algorithm}\begin{footnotesize}
\end{footnotesize}
\end{scriptsize}
\begin{remark}
Note that for every iteration $k$, Algorithm~\ref{al:IAMA for DMPC} and \ref{al:IFAMA for DMPC} only need local communication and the computations can be performed in parallel for every subsystem.
\end{remark}

%

\color{black}
We provide a lemma showing that considering Algorithm~\ref{al:IAMA for DMPC} there exists a Lipschitz constant $L(\psi)$ equal to zero. The results can be easily extended to Algorithm~\ref{al:IFAMA for DMPC}. This result is required by the proofs of the complexity upper-bounds in Corollary~\ref{co:convergence rate of inexact AMA for DMPC}, \ref{co:linear convergence rate of inexact AMA for DMPC} and \ref{co:convergence rate of inexact FAMA for DMPC}. 
\color{black}

\begin{lemma}\label{le: the lipschitz constant L(psi) is zero for IAMA for DMPC}
Let the sequence $\lambda^{k}$ be generated by Algorithm~\ref{al:IAMA for DMPC}. For all $k \geq 0$ it holds that $E^{T}\lambda^{k} = 0$ and the Lipschitz constant of the second objective in the dual problem of Problem~\ref{pr:split form for distributed MPC} $L(\psi)$ is equal to zero.
\end{lemma}

\begin{IEEEproof}
We first prove that for all $k \geq 0$, the sequence $\lambda^{k}$ satisfies $E^{T}\lambda^{k} = 0$. We know that Step~3 in Algorithm~\ref{al:IAMA for DMPC} is equivalent to the following update
\begin{align*}
\tilde{v}^{k} = \mathcal{M}\cdot \sum^{M}_{i=1} E^{T}_i \cdot \tilde{z}^{k}_i = \mathcal{M}\cdot E^{T}\cdot \tilde{z}^{k}\enspace ,
\end{align*}
with $\mathcal{M} = blkdiag(\frac{1}{|\mathcal{N}_1|}\cdot I_{1},\cdots, \frac{1}{|\mathcal{N}_i|}\cdot I_i, \cdots, \frac{1}{|\mathcal{N}_M|}\cdot I_M) = (E^{T}E)^{-1}$, where $|\mathcal{N}_i|$ denotes the number of the elements in the set $\mathcal{N}_i$, and $I_{i}$ denotes an identity matrix with the dimension of the $i$th component of $v$, denoted as $[v]_i$. From Step~5 in Algorithm~\ref{al:IAMA for DMPC}, for all $k\geq 1$ we have that
\begin{align*}
\lambda^{k} = \lambda^{k-1} + \tau(E\tilde{v}^{k}-\tilde{z}^{k})\enspace .
\end{align*}
By multiplying the matrix $E^{T}$ to both sides, we have
\begin{align*}
E^{T}\lambda^{k} 
= E^{T}\lambda^{k-1} + \tau(E^{T}E\tilde{v}^{k}-E^{T}\tilde{z}^{k})
= E^{T}\lambda^{k-1} + \tau(E^{T}E\mathcal{M} E^{T}\tilde{z}^{k}-E^{T}\tilde{z}^{k})\enspace .
\end{align*}
Since $\mathcal{M} = (E^{T}E)^{-1}$, the above equality becomes 
\begin{align*}
E^{T}\lambda^{k} = E^{T}\lambda^{k-1} + \tau(E^{T}\tilde{z}^{k}-E^{T}\tilde{z}^{k}) = E^{T}\lambda^{k-1}\enspace .
\end{align*}
From the initialization in Algorithm~\ref{al:IAMA for DMPC}, we know $E^{T}\lambda^{0} = E^{T}\cdot 0 =0$. Then by induction, we can immediately prove that for all $k \geq 0$ it holds that $E^{T}\lambda^{k} = 0$.  We can now show that for all $E^{T}\lambda = 0$, a Lipschitz constant of the second objective in the dual problem in Problem~\ref{pr:dual problem of AMA and FAMA problem} $L(\psi)$ is equal to zero. Since  $g=0$, $B = -E$ and $c=0$, then the second objective in the dual problem is equal to
 \begin{align*}
\psi(\lambda) = g^{\star}(B^{T}\lambda) - c^{T}\lambda = g^{\star}(E^{T}\lambda)=  \sup_{w} (v^{T}E^{T}\lambda - 0) = \begin{cases} 0 &\mbox{if } E^{T}\lambda = 0 \\ 
\infty & \mbox{if } E^{T}\lambda \neq 0. \end{cases} \enspace .
\end{align*}
The function $\psi(\lambda)$ is an indicator function on the nullspace of matrix $E^{T}$. For all $\lambda$ satisfying $E^{T}\lambda = 0$, the function $\psi(\lambda)$ is equal to zero. Hence, zero is a Lipschitz constant of the function $\psi(\lambda)$ for all $E^{T}\lambda = 0$. 
\end{IEEEproof}
After proving Lemma~\ref{le: the lipschitz constant L(psi) is zero for IAMA for DMPC}, we are ready to show the main theoretical properties of Algorithm~\ref{al:IAMA for DMPC} and~\ref{al:IFAMA for DMPC}.
\begin{corollary}\label{co:convergence rate of inexact AMA for DMPC}
Let $\{\lambda^k=[\lambda^{k^{T}}_1,\cdots,\lambda^{k^{T}}_M]^{T}\}$ be generated by Algorithm~\ref{al:IAMA for DMPC}. If Assumption~\ref{as: distributed problems for inexact AMA} is satisfied and the inexact solutions $\tilde{z}^{k}_i$ for all $k\geq 1$ are feasible, i.e. $\tilde{z}^{k}_i\in\mathbb{C}_i$, then for any $k\geq 1$
\begin{align}
D(\lambda^{\star}) - D\left(\frac{1}{k}\sum^{k}_{p=1} \lambda^p \right) \leq \frac{L(\nabla \phi)}{2k}\left(\|\lambda^{0}-\lambda^{\star}\|+2\sum^{k}_{p=1}\frac{\|\delta^{p}\|}{L(\nabla \phi)}\right)^{2} \enspace ,
\end{align} 
where $D(\cdot)$ is the dual function of Problem~\ref{pr:split form for distributed MPC}, $\lambda^{0}=[\lambda^{0^{T}}_1,\cdots,\lambda^{0^{T}}_M]^{T}$ and $\lambda^{\star}$ are the starting sequence and the optimal sequence of the Lagrangian multiplier, respectively, and $\delta^{p}=[\delta^{p^{T}}_1,\cdots,\delta^{p^{T}}_M]^{T}$ denotes the global error sequence. The Lipschitz constant $L(\nabla \phi)$ is equal to $\sigma_{f}^{-1}$, with $\sigma_{f}=\min \{\sigma_{f_1}, \cdots, \sigma_{f_M}\}$. 
\end{corollary}

\begin{IEEEproof}
As stated in Remark~\ref{rem:IAMA for DMPC constrains on the first objective}, Problem~\ref{pr:split form for distributed MPC} is split as follows: $f = \sum^{M}_{i=1} f_{i}(z_i)$ with the constraints $z_i\in \mathbb{C}_{i}$ for all $i=1,\cdots ,M$ and $g=0$. The matrices are $A = I$, $B = -E$ and $c = 0$. If Assumption~\ref{as: distributed problems for inexact AMA} holds, then this splitting problem satisfies Assumption~\ref{as:AMA and FAMA} with the convexity modulus $\sigma_f$. From Theorem~\ref{th:convergence rate of inexact AMA}, we know that the sequence $\{\lambda^k\}$ generated by inexact AMA in Algorithm~\ref{al:IAMA for DMPC}, satisfies the complexity upper bound in (\ref{eq:upper bound in the theorem of inexact AMA}) with $\Gamma^{k}$ and $\Lambda^{k}$ in (\ref{eq: Gamma in the bound in the theorem of inexact AMA}) and (\ref{eq: Lambda in the bound in the theorem of inexact AMA}) with $\delta^{k}=[\delta^{k^{T}}_1,\cdots,\delta^{k^{T}}_M]^{T}$ and $\theta^{k} =0$. By Lemma~\ref{le: the lipschitz constant L(psi) is zero for IAMA for DMPC}, it follows that the constant $L(\psi)$ in $\Lambda^{k}$ is equal to zero. The Lipschitz constant of the gradient of the dual objective is equal to $L(\nabla \phi)=\sigma_{f}^{-1}\cdot \rho(A) = \sigma_{f}^{-1}$ with $\sigma_{f}=\min \{\sigma_{f_1}, \cdots, \sigma_{f_M}\}$. Hence, we can simplify the complexity upper bound in (\ref{eq:upper bound in the theorem of inexact AMA}) for Algorithm~\ref{al:IAMA for DMPC} to be inequality (\ref{eq:complexity bound for inexact FAMA for DMPC}). 
%
%
\end{IEEEproof}

As we discussed in Remark~\ref{re: discussion about quadratic cost functions of DMPC for IAMA}, if the state cost functions $l_i(\cdot,\cdot)$ and $l^{f}_i(\cdot)$ in the distributed MPC problem are strictly positive quadratic functions, then the distributed optimization problem originating from the distributed MPC problem satisfies Assumption~\ref{as: AMA with quadratic f and A=I}, which according to Theorem~\ref{th:convergence rate of inexact AMA} implies a linearly decreasing upper-bound given in Corollary~\ref{co:linear convergence rate of inexact AMA for DMPC}.
\begin{corollary}\label{co:linear convergence rate of inexact AMA for DMPC}
Let $\{\lambda^k=[\lambda^{k^{T}}_1,\cdots,\lambda^{k^{T}}_M]^{T}\}$ be generated by Algorithm~\ref{al:IAMA for DMPC}. If Assumption~\ref{as: distributed problems for inexact AMA} is satisfied, the local cost function $f_i$ is a strictly positive quadratic function, and the inexact solutions $\tilde{z}^{k}_i$ for all $k\geq 1$ are feasible, i.e. $\tilde{z}^{k}_i\in\mathbb{C}_i$, then for any $k\geq 1$
\begin{align}
\|\lambda^k - \lambda^{\star}\|  \leq (1-\gamma)^{k+1}\cdot\left(\|\lambda^{0}-\lambda^{\star}\|+\sum^{k}_{p=0} (1-\gamma)^{-p-1}\cdot\frac{\|A\delta^{p}\|}{L(\nabla \phi)}\right)\enspace , 
\label{eq:upper bound in the theorem of inexact AMA with linear rate for DMPC}
\end{align} 
where $\gamma = \frac{\mathbf{\lambda}_{min}(H)}{\mathbf{\lambda}_{max}(H)}$, and $\lambda^{0}$ and $\lambda^{\star}$ are the starting sequence and the optimal sequence of the Lagrangian multiplier, respectively. The Lipschitz constant $L(\nabla \phi)$ is equal to $\sigma_{f}^{-1}$, where $\sigma_{f}=\min \{\sigma_{f_1}, \cdots, \sigma_{f_M}\}$.
\end{corollary}
\color{black}
\begin{IEEEproof}
In Algorithm~\ref{al:IFAMA for DMPC}, the variable $\hat{\lambda}^{k}_{i}$ is a linear function of $\lambda^{k}_{i}$ and $\lambda^{k-1}_{i}$. This preserves all properties shown in Lemma~\ref{le: the lipschitz constant L(psi) is zero for IAMA for DMPC} for Algorithm~\ref{al:IFAMA for DMPC}. Then, Corollary~\ref{co:linear convergence rate of inexact AMA for DMPC} can be easily proven by following the same steps as in the proof of Corollary~\ref{co:convergence rate of inexact AMA for DMPC} by replacing Theorem~\ref{th:convergence rate of inexact AMA} by Theorem~\ref{th:linear convergence rate of inexact AMA}.
\end{IEEEproof}
\color{black}

%
\begin{corollary}\label{co:convergence rate of inexact FAMA for DMPC}
Let $\{\lambda^k=[\lambda^{k^{T}}_1,\cdots,\lambda^{k^{T}}_M]^{T}\}$ be generated by Algorithm~\ref{al:IFAMA for DMPC}. If Assumption~\ref{as: distributed problems for inexact AMA} is satisfied and the inexact solutions $\tilde{z}^{k}_i$ for all $k\geq 1$ are feasible, i.e. $\tilde{z}^{k}_i\in\mathbb{C}_i$, then for any $k\geq 1$
\begin{align}\label{eq:complexity bound for inexact FAMA for DMPC}
D(\lambda^{\star}) - D(\lambda^k) \leq \frac{2L(\nabla \phi)}{(k+1)^{2}}\left(\|\lambda^{0}-\lambda^{\star}\|+2M\sum^{k}_{p=1} p\frac{\delta^{p}}{L(\nabla \phi)}\right)^{2}\enspace .
\end{align} 
where $D(\cdot)$ is the dual function of Problem~\ref{pr:split form for distributed MPC}, $\lambda^{0}$ and $\lambda^{\star}$ are the starting sequence and the optimal sequence of the Lagrangian multiplier, respectively. The Lipschitz constant $L(\nabla \phi)$ is equal to $\sigma_{f}^{-1}$, where $\sigma_{f}=\min \{\sigma_{f_1}, \cdots, \sigma_{f_M}\}$.
\end{corollary}
\color{black}
\begin{IEEEproof}
It follows from the same proof as Corollary~\ref{co:convergence rate of inexact AMA for DMPC} by replacing Theorem~\ref{th:convergence rate of inexact AMA} by Theorem~\ref{th:convergence rate of inexact FAMA}.
\end{IEEEproof}
\color{black}
\begin{remark}
For the case that all the local problems are solved exactly, i.e. $\delta^{k}_{i}=0$, Algorithm~\ref{al:IAMA for DMPC} and Algorithm~\ref{al:IFAMA for DMPC} reduce to standard AMA and FAMA, and converge to the optimal solution at the rate of the complexity upper-bounds.
\end{remark}
\begin{remark}\label{re: sufficient condition on the local errors for IAMA for DMPC}
The sufficient conditions on the errors for convergence given in Corollary~\ref{co:sufficient conditions on the errors for IAMA}, \ref{co:sufficient conditions on the errors for IAMA with linear convergence rate} and \ref{co:sufficient conditions on the errors for IFAMA} can be directly extended to the error sequence $\{\delta^{k}\}$.
%
%
\end{remark}

\color{black}
\subsection{Certification of the number of local iterations for convergence}\label{se: Certificate of the number of iterations for local problems}


We have shown that the inexact distributed optimization algorithms in Algorithm~\ref{al:IAMA for DMPC} and \ref{al:IFAMA for DMPC} allow one to solve the local problems, i.e. Step~1 in Algorithm~\ref{al:IAMA for DMPC} and \ref{al:IFAMA for DMPC}, inexactly. In this section, we will address two questions: which algorithms are suitable for solving the local problems; and what termination conditions for the local algorithms guarantee that the computational error of the local solution satisfies the sufficient conditions on the errors for the global distributed optimization algorithms. 

 We apply the proximal gradient method for solving the local problems in Step~1 in Algorithm~\ref{al:IAMA for DMPC} and \ref{al:IFAMA for DMPC}, and propose an approach to certify the number of iterations for their solution, by employing a warm-start strategy and the complexity upper-bounds of the proximal gradient method. The approach guarantees that the local computational errors $\delta^{k}_{i}$ decrease with a given rate, that satisfies the sufficient conditions derived from Corollary~\ref{co:convergence rate of inexact AMA for DMPC}, \ref{co:linear convergence rate of inexact AMA for DMPC} and \ref{co:convergence rate of inexact FAMA for DMPC}, ensuring convergence of the inexact distributed optimization algorithm to the optimal solution. We define a decrease function $\alpha^{k}$ satisfying the sufficient conditions, for example $\alpha^{k} = \alpha^{0} \cdot \frac{1}{k^{2}}$, where $\alpha^{0}$ is a positive number.
 
\subsubsection{Gradient method}

 The local problems in Step~1 in Algorithm~\ref{al:IAMA for DMPC} and \ref{al:IFAMA for DMPC} are optimization problems with strongly convex cost functions and convex constraints. From Corollary~\ref{co:convergence rate of inexact AMA for DMPC}, \ref{co:linear convergence rate of inexact AMA for DMPC} and \ref{co:convergence rate of inexact FAMA for DMPC}, we know that the inexact solution $\tilde{z}^{k}_{i}$ needs to be a feasible solution subject to the local constraint $\mathbb{C}_{i}$, i.e., $\tilde{z}^{k}_{i} \in \mathbb{C}_{i}$ for all $k>0$. Therefore, a good candidate algorithm for solving the local problems should have the following three properties: the algorithm can solve convex optimization problems efficiently; if the algorithm is stopped early, i.e., only a few number of iterations are implemented, the sub-optimal solution is feasible with respect  to the local constraint $\mathbb{C}_{i}$; and there exists a certificate on the number of iterations to achieve a given accuracy of the sub-optimal solution. Gradient methods satisfy these requirements, have simple and efficient implementations, and offer complexity upper-bounds on the number of iterations \cite{beck_fISTA_2009}. These methods have been studied in the context of MPC in \cite{richter_computational_2012}, \cite{giselsson_accelerated_2013} and \cite{pu_FAMAMPC_2014}. 
 
 
  We apply the proximal gradient method in Algorithm~\ref{al:fast gradient method for inner-loop} for solving the local problems in Step~1 in Algorithm~\ref{al:IAMA for DMPC} and \ref{al:IFAMA for DMPC}. The local optimization problems at iteration $k$ are parametric optimization problems with the parameter $\lambda^{k-1}_i$. We denote the optimal function as 
\begin{equation}\label{eq: the optimal function for the parametric local problem}
\textbf{z}^{\star}_{i}(\lambda_i) := \mbox{argmin}_{z_i\in \mathbb{C}_{i}}\lbrace  f_i(z_i) + \langle \lambda_{i}, -z_i\rangle\rbrace \enspace .
\end{equation} 
The solution of the optimal function at $\lambda^{k-1}_i$ is denoted as $z^{k, \star}_{i} := \textbf{z}^{\star}_{i}(\lambda^{k-1}_i)$. The function $\textbf{z}^{\star}_{i}(\cdot)$ has a Lipschitz constant $L(\textbf{z}^{\star}_{i})$ satisfying as $\|\textbf{z}^{\star}_{i}(\lambda_{i_1}) - \textbf{z}^{\star}_{i}(\lambda_{i_1})\|\leq L(\textbf{z}^{\star}_{i})\cdot\|\lambda_{i_1} - \lambda_{i_2}\|$ for any $\lambda_{i_1}$ and $\lambda_{i_2}$. Motivated by the fact that the difference between the parameters $\lambda^{k-1}_i$ and $\lambda^{k}_i$ is limited and measurable for each $k$, i.e. $\beta^{k}_i = \|\lambda^{k-1}_i-\lambda^{k}_i\| = \tau(E_i\tilde{v}^{k-1}-\tilde{z}^{k-1}_i)$, we use a warm-starting strategy to initialize the local problems, i.e. we use the solution $\tilde{z}^{k-1}_{i}$ from the previous step $k-1$ as the initial solution for Algorithm~\ref{al:fast gradient method for inner-loop} for step $k$. 

\begin{scriptsize}
\begin{algorithm}
\caption{Gradient method for solving Step~1 in Algorithm~\ref{al:IAMA for DMPC} at iteration $k$}
\begin{algorithmic}
\REQUIRE Initialize $\alpha^{k}=\alpha^{0} \cdot \frac{1}{k^{2}}$, $\beta^{k} = \|\tau(E_i\tilde{v}^{k-1}-\tilde{z}^{k-1}_i)\|$, $\lambda^{k-1}_{i}$, $z^{k,0}_{i} = \tilde{z}^{k-1}_{i}$ and $\tau_i < \frac{1}{L(\nabla f_i)}$
\STATE Compute $J_k$ satisfying (\ref{eq: condition for number of iterations for solving the local problems Jk})
\FOR {$j=1,2,\cdots, J_k$}
\STATE $z^{k, j}_{i} =  \mbox{Proj}_{\mathbb{C}_i}(z^{k, j-1}_{i} - \tau (\nabla f_i(z^{k, j-1}_{i}) -\lambda^{k-1}_{i}))$
\ENDFOR
\STATE $\tilde{z}^{k}_{i} \gets z^{k, J_k}_{i}$
\end{algorithmic}
\label{al:fast gradient method for inner-loop}
\end{algorithm}\begin{footnotesize}
\end{footnotesize}
\end{scriptsize}

Note that we initialize the vectors $\tilde{v}^{k-1}$, $\tilde{z}^{k-1}_i$ and $\tilde{z}^{k-1}_{i}$ for $k=1$ in Algorithm~\ref{al:IAMA for DMPC} to be zero vectors.

\begin{proposition}[Proposition 3 in \cite{schmidt_convergence_2011}]\label{pr:convergence rate of ISTA}
Let $z^{k, j}_{i}$ be generated by Algorithm~\ref{al:fast gradient method for inner-loop}. If Assumption~\ref{as: distributed problems for inexact AMA} holds, then for any $j\geq 0$ we have:
\begin{align}
\|z^{k, j}_{i}-z^{k, \star}_{i}\|\leq \|z^{k,0}_{i}-z^{k, \star}_{i}\|\cdot(1-\gamma)^{j}\enspace ,
\end{align}
where $\gamma = \frac{\sigma_{f_i}}{L(\nabla f_i)}$ and $z^{k,0}_{i}$ and $z^{k, \star}_{i}$ denote the initial sequence of  Algorithm~\ref{al:inexact ISTA} and the optimal solution of the problem in Step~6 in Algorithm~\ref{al:IAMA for DMPC} at iteration $k$, respectively.
\end{proposition}

\subsubsection{Termination condition on the number of iterations for solving local problems}\label{se: Termination condition on the number of iterations for solving local problems}

Methods for bounding the number of iterations to reach a given accuracy have been studied e.g. in \cite{giselsson_execution_2012}, \cite{patrinos_accelerated_2014} and \cite{dinh_fast_2012}. In \cite{giselsson_execution_2012} and  \cite{patrinos_accelerated_2014}, the authors proposed dual decomposition based optimization methods for solving quadratic programming problems and presented termination conditions to guarantee a prespecified accuracy. However, these methods do not directly guarantee feasibility of the sub-optimal solution. One approach is to tighten constraints to ensure feasibility, which can be conservative in practice. In \cite{dinh_fast_2012}, the authors propose an inexact decomposition algorithm for solving distributed optimization problems by employing smoothing techniques and an excessive gap condition as the termination condition on the number of iterations to achieve a given accuracy. To certify the termination condition, this method requires to measure the values of the global primal and dual functions on-line, which requires full communication on the network and is not satisfied in our distributed framework. In addition, this method does not provide any algorithms for solving the local problems. 

By employing the complexity upper-bounds in Proposition~\ref{pr:convergence rate of ISTA} for Algorithm~\ref{al:fast gradient method for inner-loop}, we propose a termination condition in (\ref{eq: condition for number of iterations for solving the local problems Jk}) to find the number of iterations $J_k$, which guarantees that the local computational error is upper-bounded by the predefined decrease function $\alpha^{k}$, i.e. $\|\delta^{k}_{i}\| \leq \alpha^{k}$.
\begin{lemma}\label{le: the local computation errors satisfy a given decreasing rate}
If the number of iterations $J_k$ in Algorithm~\ref{al:fast gradient method for inner-loop} satisfies 
\begin{equation}\label{eq: condition for number of iterations for solving the local problems Jk}
J_k \geq \lceil\log_{(1-\gamma)} \frac{\alpha^{k}}{\alpha^{k-1}+L(\textbf{z}^{\star}_i)\beta^{k}}\rceil
\end{equation}
for all $k\geq 1$, then the computational error for solving the local problem in Step~6 in Algorithm~\ref{al:IAMA for DMPC} $\delta^{k}_{i}$ satisfies $\|\delta^{k}_{i}\| \leq \alpha^{k}$.
\end{lemma}
\begin{IEEEproof}
We will prove Lemma~\ref{le: the local computation errors satisfy a given decreasing rate} by induction.
\begin{itemize}
\item Base case: For $k = 1$, the vectors $\tilde{v}^{k-1}$, $\tilde{z}^{k-1}_i$ and $\tilde{z}^{k-1}_{i}$ are initialized as zero vectors. By Proposition~\ref{pr:convergence rate of inexact ISTA} and the fact $z^{1, 0}_{i} = \tilde{z}^{0}_{i}=0$, we know
\begin{align*}
\|z^{1, J_1}_{i}-z^{1, \star}_{i}\|
\leq \|z^{1, 0}_{i}-z^{1, \star}_{i}\|\cdot(1-\gamma)^{J_1}= \|0 - z^{1, \star}_{i}\|\cdot(1-\gamma)^{J_1}\enspace .
\end{align*}
Due to the definition of the function $\alpha^{k}$, it follows that the term above is upper-bounded by $\alpha^{0}\cdot(1-\gamma)^{J_1}$. Using the fact that $\beta^{0} = \|\tau(E_i\tilde{v}^{0}-\tilde{z}^{0}_i)\|=0$ and $J_1$ satisfies (\ref{eq: condition for number of iterations for solving the local problems Jk}), it is further upper-bounded by $\alpha^{1}$;
\begin{align*}
\|\delta^{1}\| = \|\tilde{z}^{1}_{i}-z^{1, \star}_{i}\|= \|z^{1, J_1}_{i}-z^{1, \star}_{i}\|\leq \alpha^{1}\enspace .
\end{align*}
\item Induction step: Let $l\geq 1$ be given and suppose that $\|\delta^{l}\| \leq \alpha^{l}$. We will prove that $\|\delta^{l+1}\| \leq \alpha^{l+1}$
By Proposition~\ref{pr:convergence rate of inexact ISTA} and the warm-starting strategy, i.e. $z^{k, 0}_{i} = \tilde{z}^{k-1}_{i}=z^{k-1, J_{k-1}}_{i}$, we know
\begin{align*}
\|\delta^{l+1}\| &= \|z^{l+1, J_{g+1}}_{i}-z^{l+1, \star}_{i}\|\\
&\leq \|z^{l+1, 0}_{i}-z^{l+1, \star}_{i}\|\cdot(1-\gamma)^{J_{l+1}}\\
&= \|z^{g, J_{l}}_{i} - z^{l+1, \star}_{i}\|\cdot(1-\gamma)^{J_{l+1}}\\
&\leq (\|z^{l, J_{l}}_{i}-z^{l, \star}_{i}\| + \| z^{l, \star}_{i} - z^{l+1, \star}_{i}\|)\cdot(1-\gamma)^{J_{l+1}}\\
&\leq (\delta^{l} + L(\textbf{z}^{\star}_i)\cdot \beta^{l+1})\cdot(1-\gamma)^{J_{l+1}}\enspace .
\end{align*}
Due to the induction assumption and the fact that $J_l$ satisfies (\ref{eq: condition for number of iterations for solving the local problems Jk}), it follows that $\|\delta^{l+1}\| \leq \alpha^{l+1}$.
\end{itemize}
We conclude that by the principle of induction, it holds that $\|\delta^{k}\| \leq \alpha^{k}$ for all $k\geq 1$.
\end{IEEEproof}

\begin{corollary}\label{co:convergence of inexact AMA for DMPC with on-line certification}
 If Assumption~\ref{as: distributed problems for inexact AMA} holds and the decrease rate of the function $\alpha^{k}$ satisfies the corresponding sufficient conditions presented in Corollary~\ref{co:sufficient conditions on the errors for IAMA} and \ref{co:sufficient conditions on the errors for IFAMA}, then Algorithm~\ref{al:IAMA for DMPC} and \ref{al:IFAMA for DMPC}  converge to the optimal solution, with Algorithm~\ref{al:fast gradient method for inner-loop} solving the local problem in Step~1. Furthermore, if the local cost function $f_i$ is a strictly positive quadratic function, and the decrease rate of the function $\alpha^{k}$ satisfies the sufficient conditions presented in Corollary~\ref{co:sufficient conditions on the errors for IAMA with linear convergence rate}, then Algorithm~\ref{al:IAMA for DMPC} converges to the optimal solution, with Algorithm~\ref{al:fast gradient method for inner-loop} solving the local problem in Step~1. 
\end{corollary}
%
%
\begin{remark}
All the information required by the proposed on-line certification method, i.e.,by Algorithm~\ref{al:fast gradient method for inner-loop}, as well as the condition for $J_k$ in (\ref{eq: condition for number of iterations for solving the local problems Jk}), can be obtained on-line and locally.
\end{remark}

\subsubsection{Computation of the Lipschitz constant $L(\textbf{z}^{\star}_i)$}

In the above proposed on-line certification method, the Lipschitz constant of the optimal solution function $\textbf{z}^{\star}_i(\lambda^{k-1}_{i})$, $L(\textbf{z}^{\star}_i)$, plays an important role. While it is generally difficult to compute this Lipschitz constant, it can be computed for special cases, such as positive quadratic functions.
\begin{lemma}\label{le: the Lipschitz constant of zstar}
Let the local cost function be a quadratic function, i.e. $f_i(z_i) = \frac{1}{2}z^{T}_{i}H_iz_i + h^{T}_iz_i$ with $H_i\succ 0$. A Lipschitz constant of the function $\textbf{z}^{\star}_{i}(\lambda_i)$ defined in (\ref{eq: the optimal function for the parametric local problem}) is given by $\frac{1}{\rho_{min}(H_i)}$, i.e. 
\begin{equation}
\|\textbf{z}^{\star}_{i}(\lambda_{i_1}) -\textbf{z}^{\star}_{i}(\lambda_{i_2})\| \leq \frac{1}{\rho_{min}(H_i)}\cdot \|\lambda_{i_1}-\lambda_{i_2}\|.
\end{equation} 
\begin{IEEEproof}
Since $H_i\succ 0$, we can define $H_i = D\cdot D^{T}$ with $D$ invertible, which implies
\begin{align*}
\textbf{z}^{\star}_{i}(\lambda_{i}) 
& =  \operatorname{argmin}_{z_i\in \mathbb{C}_{i}} \;\frac{1}{2}z^{T}_{i}H_iz_i + (h_i-\lambda_{i})^{T}z_i\\
& = \operatorname{argmin}_{z_i\in \mathbb{C}_{i}} \;\frac{1}{2} \|D^{T}z_i + D^{-1}(h_i-\lambda_{i})\|^{2} \\
\end{align*}
Let $v = D^{T}z_i$. The optimization problem above becomes 
\begin{align*}
\textbf{v}^{\star}(\lambda_{i}) 
& =  \operatorname{argmin}_{D^{T^{-1}}v\in \mathbb{C}_{i}} \; \;\frac{1}{2} \|v + D^{-1}(h_i-\lambda_{i})\|^{2} \enspace ,
\end{align*}
which can be seen as the projection of the point $D^{-1}(h_i-\lambda_{i})$ onto the set $\bar{\mathbb{C}}_{i}:=\{v\mid D^{-1}v_i\in \mathbb{C}_{i}\}$. Since $\mathbb{C}_{i}$ is convex, then $\bar{\mathbb{C}}_{i}$ is convex as well. It follows directly from Proposition~2.2.1 in \cite{Bertsekas_Convex_2003} that
\begin{align*}
\|\textbf{v}^{\star}(\lambda_{i_1}) -\textbf{v}^{\star}(\lambda_{i_2})\| \leq \|D^{-1}\cdot (\lambda_{i_1}-\lambda_{i_2})\| \enspace .
\end{align*}
By $z_i = D^{T^{-1}}v$, we get
\begin{align*}
\|\textbf{z}^{\star}_i(\lambda_{i_1}) -\textbf{z}^{\star}_i(\lambda_{i_2})\| &\leq \|D^{-1}\|\cdot\|D^{-1}\cdot (\lambda_{i_1}-\lambda_{i_2})\| \\
& \leq \|D^{-1}\|^{2}\cdot\|\lambda_{i_1}-\lambda_{i_2}\| \\
&\leq \frac{1}{\rho_{min}(H_i)}\cdot \|\lambda_{i_1}-\lambda_{i_2}\|\enspace .
\end{align*}
\end{IEEEproof}
\end{lemma}
\color{black}

\section{Numerical example}\label{SEC:example}
This section illustrates the theoretical findings of the paper and demonstrates the performance of inexact AMA by solving a randomly generated distributed MPC problem with $40$ sub-systems. For this example, we assume that the sub-systems are coupled only in the control input:
%
\begin{equation*}
x_{i}(t+1) = A_{ii}x_{j}(t) + \sum_{j\in \mathcal{N}_{i}} B_{ij}u_{j}(t) \quad i=1,2,\cdots ,M,
\end{equation*}
The input-coupled dynamics allow us to eliminate the states of the distributed MPC problem, such that the optimization variable in the distributed optimization problems is the control sequence $u=[u^{T}_1,\cdots,u^{T}_M]^{T}$, with $u_{i}=[u^{T}_{i}(0),u^{T}_{i}(1),\cdots,u^{T}_{i}(N)]^{T}$. Examples with this structure include systems sharing one resource, e.g. a water-tank system or an energy storage system. 

We randomly generate a connected network with $40$ agents. Each sub-system has three states and two inputs. The dynamical matrices $A_{ii}$ and $B_{ij}$ are randomly generated, i.e. generally dense, and the local systems are controllable. The input constraint $\mathbb{U}_i$ for sub-system $i$ is set to be $\mathbb{U}_i = \{u_i(t)| -0.4\leq u_i(t)\leq 0.3\}$. The horizon of the MPC problem is set to be $N=11$. The local cost functions are set to be quadratic functions, i.e. $l_i(x_i(t),u_i(t)) = x^{T}_{i}(t)Qx_i(t)+u^{T}_{i}(t)Ru_i(t)$ and $l^{f}_i(x_i(N)) = x^{T}_{i}(N)Px_i(N)$, where $Q$, $R$ and $P$ are identity matrices. Therefore, the distributed optimization problem resulting from the distributed MPC satisfies Assumption~\ref{as: distributed problems for inexact AMA}, the local cost functions $f_i$ are strictly positive quadratic functions, and the results in Corollary~\ref{co:linear convergence rate of inexact AMA for DMPC} hold. The initial states $\bar{x}_{i}$ are chosen, such that more than $70\%$ of the elements of the vector $u^{\star}$ are at the constraints. 

\begin{figure}[htbp]
   \centering
      \includegraphics[width=0.5\linewidth]{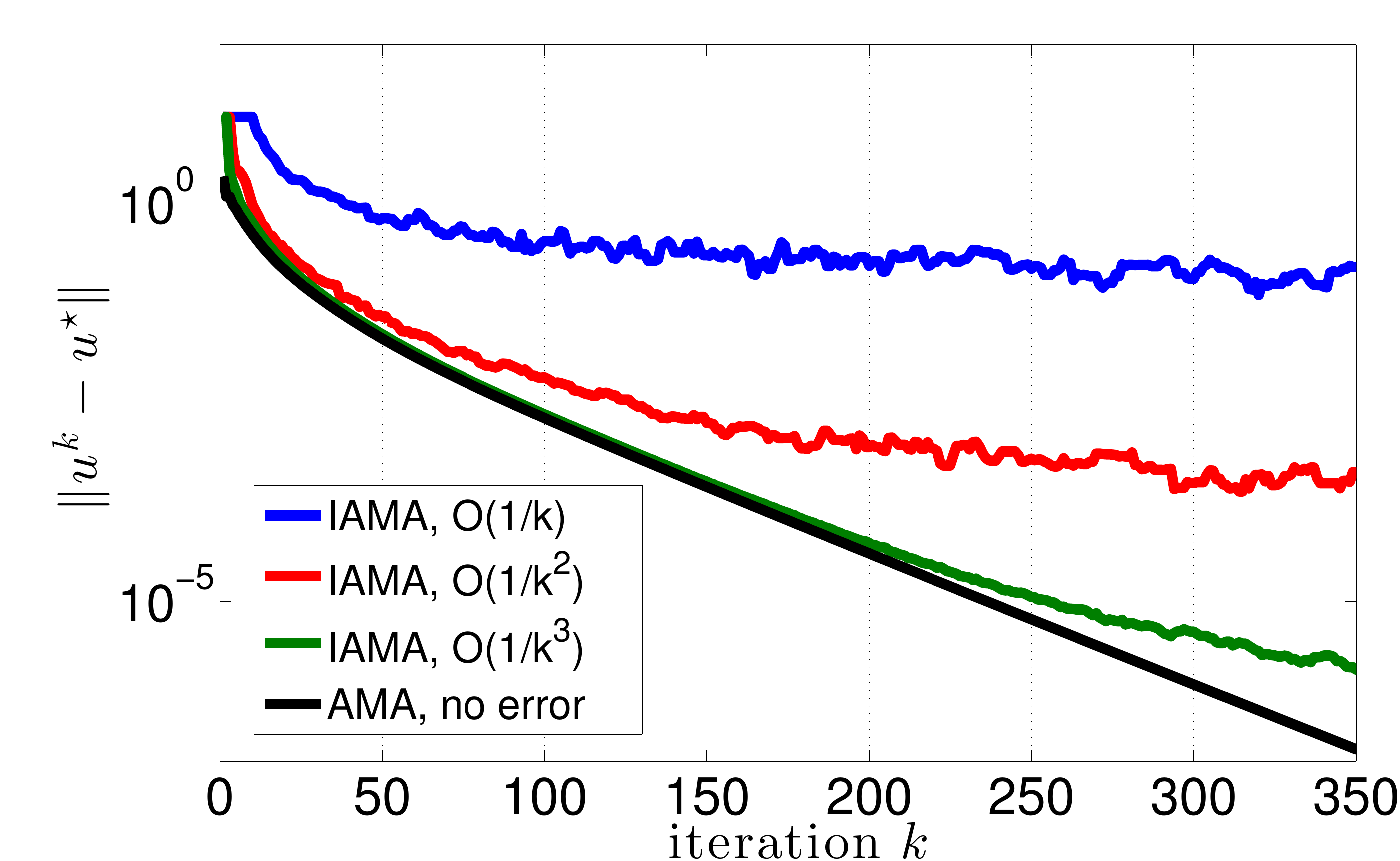}
   \caption{Comparison of the performance of AMA and inexact AMA (IAMA) with the errors decreasing at pre-defined rates.}
   \label{fig: the comparison of AMA and inexact AMA}
\end{figure}

In Fig.~\ref{fig: the comparison of AMA and inexact AMA}, we demonstrate the convergence performance of inexact AMA for solving the distributed optimization problem in Problem~\ref{pr:split form for distributed MPC}, originating from the randomly generated distributed MPC problem, applying Algorithm~\ref{al:IAMA for DMPC}. In this simulation, we compare the performance of inexact AMA with three different kinds of errors for $\delta^{k}$ with exact AMA, for which the errors are equal to zero. Note that these errors are synthetically constructed to specify different error properties. We solve the local problems to optimality and then add errors with predefined decreasing rates to the local optimal solution, ensuring that the solution remains primal feasible. The black line shows the performance of exact AMA. The blue, red and green lines show the performance of inexact AMA, where the errors $\delta^{k}$ are set to be decreasing at the rates of $O(\frac{1}{k})$, $O(\frac{1}{k^{2}})$ and $O(\frac{1}{k^{3}})$, respectively. Note that all three errors satisfy the sufficient condition for convergence in Corollary~\ref{co:sufficient conditions on the errors for IAMA with linear convergence rate}. We can observe that as the number of iterations $k$ increases, the differences $\|u^{k}-u^{\star}\|$ decrease for all the cases, however, the convergence speed is quite different. For the exact AMA algorithm (black line), it decreases linearly, which supports the results in Corollary~\ref{co:linear convergence rate of inexact AMA for DMPC}. For the three cases for inexact AMA (blue, red and green lines), we can see that the differences $\|u^{k}-u^{\star}\|$ decrease more slowly than for exact AMA, and the decrease rates correspond to the decrease rate of the errors, which supports the theoretical findings in Corollary~\ref{co:sufficient conditions on the errors for IAMA with linear convergence rate}.

The second simulation illustrates the convergence properties of inexact AMA, where the proximal gradient method in Algorithm~\ref{al:fast gradient method for inner-loop} is applied to solve the local problems in Step~2 in Algorithm~\ref{al:IAMA for DMPC}. In this experiment Algorithm~\ref{al:fast gradient method for inner-loop} is stopped after the number of iterations providing that the local computation error $\delta^{k}_i$ decreases at a certain rate. The error decrease rate is selected to be $O(\frac{1}{k})$, i.e., the decrease function $\alpha^{k}$ is set to be $\alpha^{k} = \alpha^{0} \cdot \frac{1}{k}$ and thereby satisfies the second sufficient condition in Corollary~\ref{co:sufficient conditions on the errors for IAMA with linear convergence rate}. In order to ensure $\|\delta^{k}_i\| \leq \alpha^{k}$, the number of iterations for the proximal gradient method $J_K$ in Algorithm~\ref{al:fast gradient method for inner-loop} is chosen according to the certification method presented in Section~\ref{se: Certificate of the number of iterations for local problems} such that condition  (\ref{eq: condition for number of iterations for solving the local problems Jk}) is satisfied. Note that we use a warm-starting strategy for the initialization of Algorithm~\ref{al:fast gradient method for inner-loop}. 
 
\begin{figure}[htbp]
   \centering
      \includegraphics[width=0.5\linewidth]{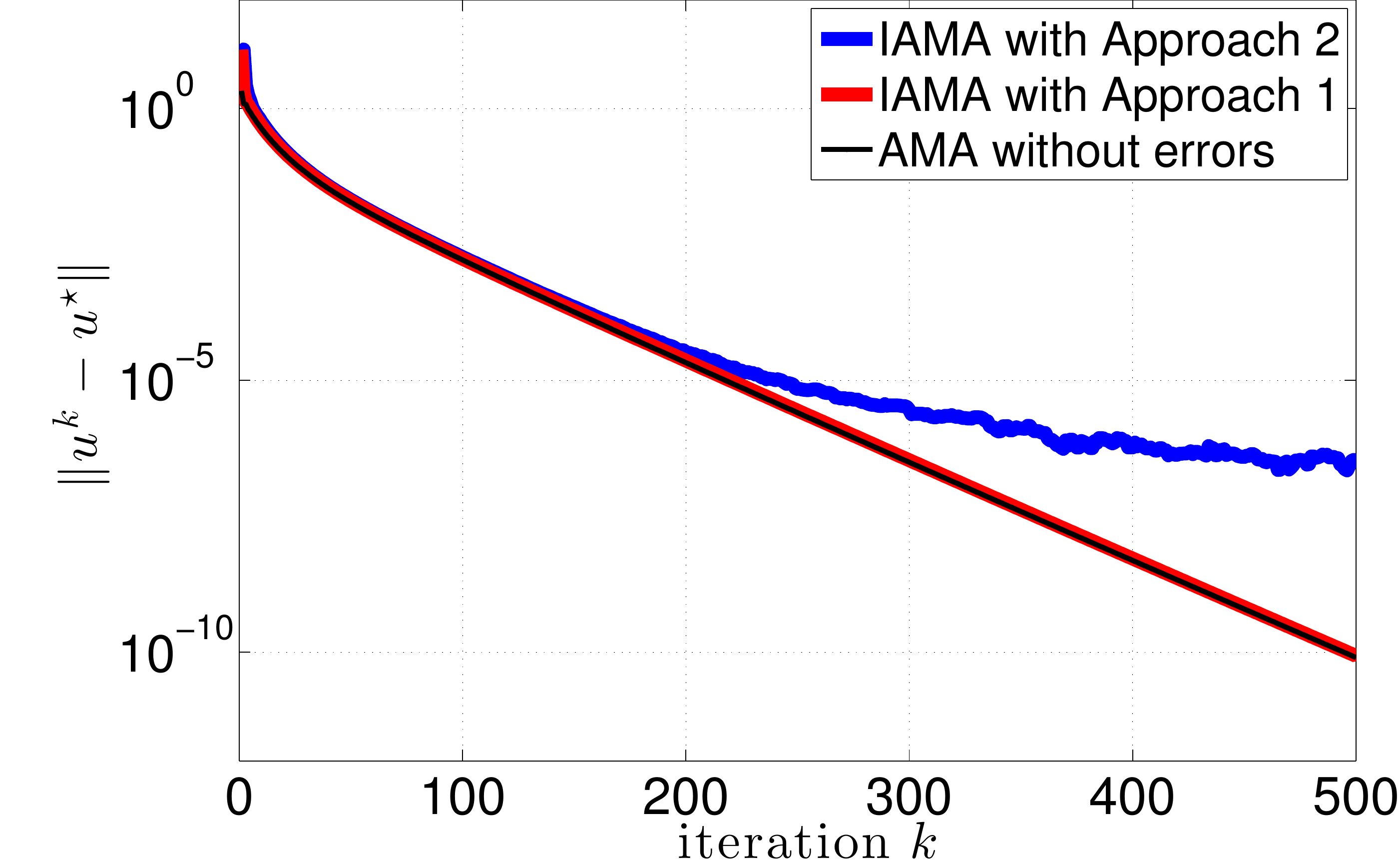}
   \caption{Comparison of the performance of AMA and inexact AMA (IAMA) with the proximal-gradient method to solve local problems, where the number of iterations is chosen according to two approaches:  \textbf{Approach~1} uses a bound on the minimal number, i.e. the termination condition proposed in (\ref{eq: condition for number of iterations for solving the local problems Jk}); and \textbf{Approach~2} computes the exact minimal number, which requires the optimal solution of the local problem $z^{k,\star}_{i}$ at each iteration.}  
   \label{fig: the comparison inexact AMA with PGM}
\end{figure}

 Fig.~\ref{fig: the comparison inexact AMA with PGM} shows the comparison of the performance of exact AMA and inexact AMA. We can observe that the black (exact AMA) and red lines basically overlap (inexact AMA with Algorithm~\ref{al:fast gradient method for inner-loop} solving local problems with the numbers of iterations $J_k$ satisfying (\ref{eq: condition for number of iterations for solving the local problems Jk})). Inexact AMA converges to the optimal solution as the iterations increase, and shows almost the same performance as exact AMA.

\begin{figure}[htbp]
   \centering
      \includegraphics[width=0.5\linewidth]{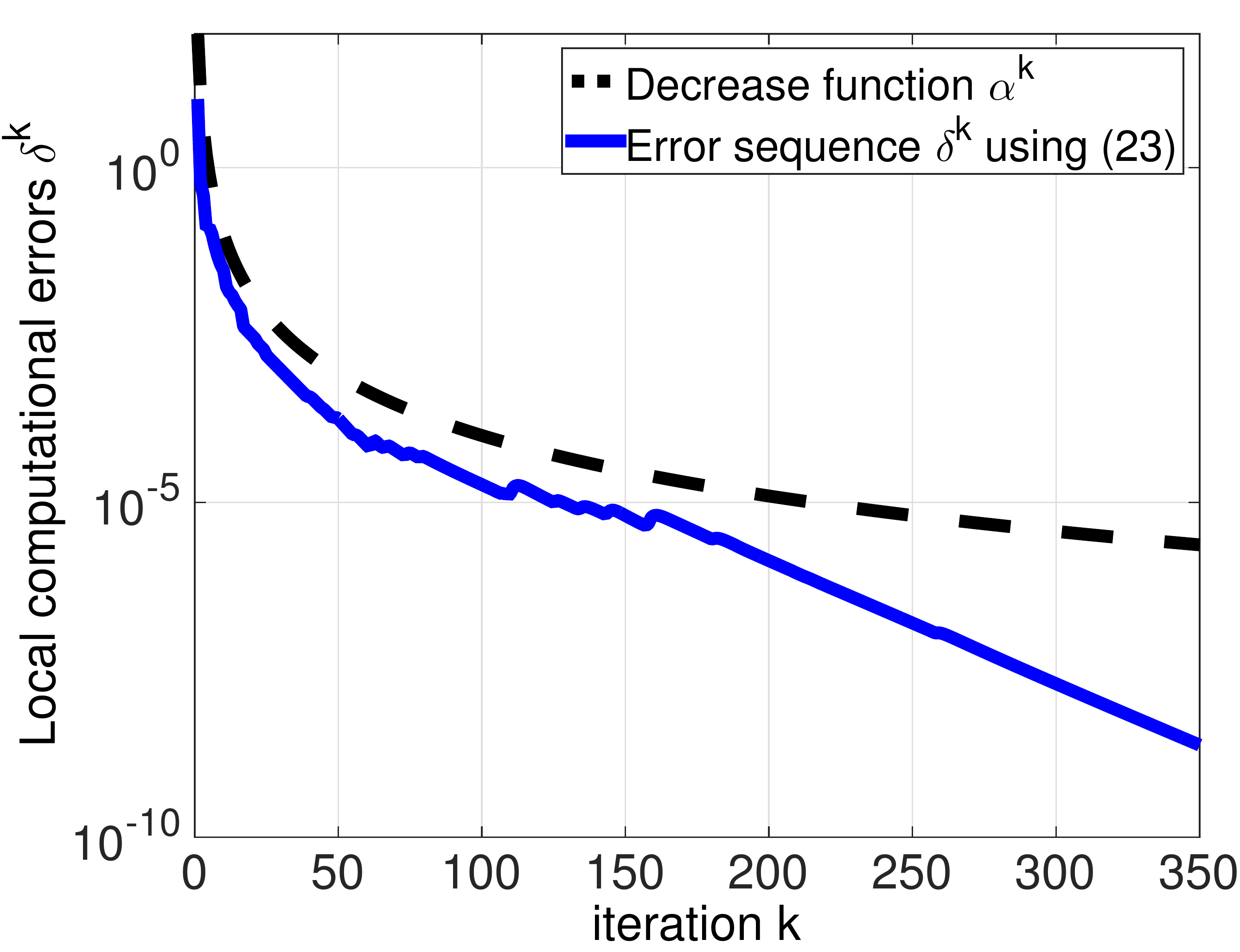}
   \caption{Error sequence $\delta^{k}$ in inexact AMA using the proximal-gradient method for solving the local problems with the numbers of iterations satisfying (\ref{eq: condition for number of iterations for solving the local problems Jk}).}
   \label{fig: the inner error of IAMA with PGM}
\end{figure}

Fig.~\ref{fig: the inner error of IAMA with PGM} shows the corresponding local error sequence $\delta^{k}_i$, where the number of iterations $J_k$ for Algorithm~\ref{al:fast gradient method for inner-loop} satisfies the condition in (\ref{eq: condition for number of iterations for solving the local problems Jk}). We can observe that the global error sequence $\delta^{k}=[\delta^{k}_1,\cdots , \delta^{k}_M]$ is upper-bounded by the decease function $\alpha^{k}$. As $k$ is small, the upper-bound $\alpha^{k}$ is tight to the error sequence. As $k$ increases, the error decreases faster and the bound becomes loose. 

\color{black}

\begin{figure}[htbp]
   \centering
      \includegraphics[width=0.6\linewidth]{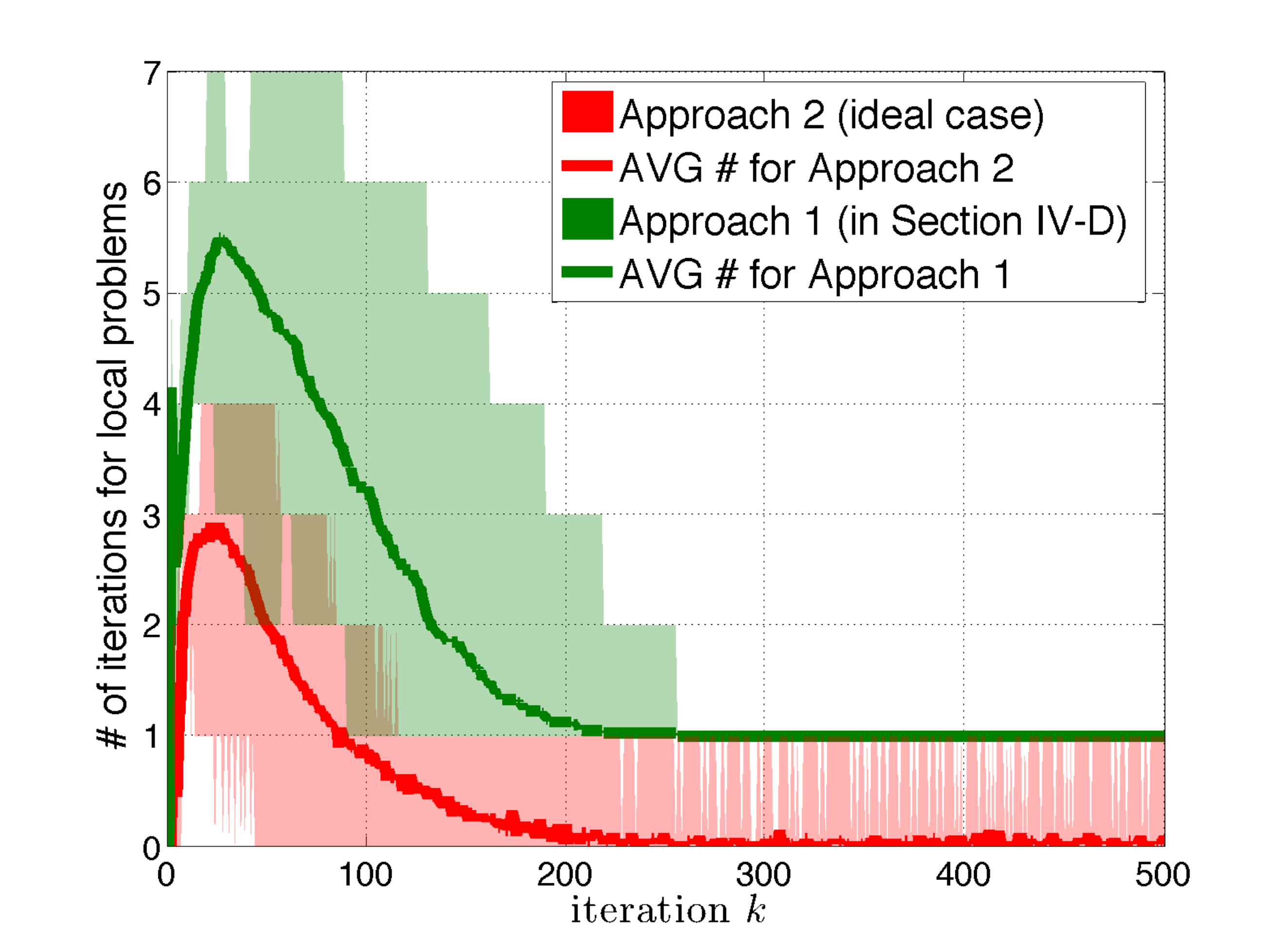}
   \caption{Comparison of the numbers of iterations for Algorithm~\ref{al:fast gradient method for inner-loop}, using two approaches:  \textbf{Approach~1} uses a bound on the minimal number, i.e. the termination condition proposed in (\ref{eq: condition for number of iterations for solving the local problems Jk}); and \textbf{Approach~2} computes the exact minimal number, which requires the optimal solution of the local problem $z^{k,\star}_{i}$ at each iteration.} 
   \label{fig: the number of iter for local problem of IAMA with PGM}
\end{figure}

Fig.~\ref{fig: the number of iter for local problem of IAMA with PGM} shows the comparison of the numbers of iterations for Algorithm~\ref{al:fast gradient method for inner-loop}, computed using two different approaches. \textbf{Approach~1} uses the termination condition proposed in Section~\ref{se: Certificate of the number of iterations for local problems}. In \textbf{Approach~2}, we first compute the optimal solution of the local problem $z^{k,\star}_{i}$ and then run the proximal gradient method to find the smallest number of iterations providing that the difference of the local sub-optimal solution satisfies the decrease function $\alpha^{k}$, i.e. $\|z^{k,j}_{i} - z^{k,\star}_{i}\|\leq \alpha^{k}$. \textbf{Approach~2} is therefore the exact minimal number, whereas \textbf{Approach~1} uses a bound on the minimal number. Note that the second approach guarantees $\|\delta^{k}_i\| \leq \alpha^{k}$ for all $k$, however, this method  is not practically applicable, since the optimal solution $z^{k,\star}_{i}$ is unknown. Its purpose is merely to compare with the proposed certification method and to show how tight the theoretical bound in (\ref{eq: condition for number of iterations for solving the local problems Jk}) is. For both techniques, we use a warm-starting strategy for initialization of the proximal gradient method to solve the local problems for each $k$ in Algorithm~\ref{al:IAMA for DMPC}. In Fig.~\ref{fig: the number of iter for local problem of IAMA with PGM}, the green line and region result from the termination condition proposed in Section~\ref{se: Certificate of the number of iterations for local problems}, and the pink line and region result from the second approach. The solid green and red lines show the average value of the numbers of iterations for the proximal gradient method for solving the local problems over the $40$ sub-systems. The upper and lower boundaries of the regions show the maximal and minimal number of iterations, respectively. The maximal number of iterations for the proposed certification method(green region) is equal to $7$, while for the second method (the red region) it is equal to $4$. Fig.~\ref{fig: the number of iter for local problem of IAMA with PGM} shows that the certification approach in (\ref{eq: condition for number of iterations for solving the local problems Jk}), which can be performed locally, is reasonably tight and the provided number of iterations is close to the minimal number of iterations required to satisfy the desired error. 

\color{black}

\section{Appendix}\label{SEC: appendix}

\subsection{Proof of Lemma~\ref{le:the equivalence between Inexact AMA and inexact ISTA}}\label{SEC: appendix proof of le:the equivalence between Inexact AMA and inexact ISTA}

\begin{IEEEproof}
In order to show the equivalence, we prove that Step~1, 2 and 3 in Algorithm~\ref{al:Inexact AMA} are equivalent to Step 1 in Algorithm~\ref{al:inexact ISTA}, i.e. the following equality holds:
\begin{equation}\label{eq: equivalence of inexact AMA and inexact ISTA on dual}
\lambda^{k} = \mbox{prox}_{\tau \psi, \epsilon^{k}}(\lambda^{k-1} - \tau (\nabla \phi(\lambda^{k-1}) + e^{k}))
\end{equation}
with $e^{k} = A\delta^{k}$ and $\epsilon^{k}=\tau ^{2}L(\psi)\|B\theta^{k}\| + \frac{\tau^{2}}{2}\|B\theta^{k}\|^{2}$. Step~2 in Algorithm~\ref{al:Inexact AMA} implies:
\begin{equation*}
B^{T}\lambda^{k-1} + \tau B^{T}(c-A\tilde{x}^{k}-Bz^{k}) \in \partial g(z^{k}),
\end{equation*}
where $z^{k}= \mbox{argmin}_{z} \;\lbrace g(z) + \langle \lambda^{k-1},-Bz\rangle+ \frac{\tau}{2}\| c-A\tilde{x}^{k+1}$
$ - Bz\|^{2}\rbrace=\tilde{z}^{k}-\theta^{k}$. From the property of the conjugate function $p\in \partial f(q) \Leftrightarrow q\in \partial f^{\star}(p)$, it follows: 
\begin{equation*}
z^{k} \in \partial g^{\star}(B^{T}\lambda^{k-1} + \tau B^{T}(c-A\tilde{x}^{k}-Bz^{k})).
\end{equation*}
By multiplying with $B$ and subtracting $c$ on both sides, we obtain:
\begin{equation*}
Bz^{k} -c \in B\partial g^{\star}(B^{T}\lambda^{k-1} + \tau B^{T}(c-A\tilde{x}^{k}-Bz^{k})) -c.
\end{equation*}
By multiplying with $\tau$ and adding $\lambda^{k-1} + \tau (c-A\tilde{x}^{k}-Bz^{k})$ on both sides, we get:
\begin{align*}
\lambda^{k-1} -\tau A\tilde{x}^{k} \in  \; \tau B\partial g^{\star}(B^{T}\lambda^{k-1} + \tau B^{T}(c-A\tilde{x}^{k}-Bz^{k})) -\tau c+\lambda^{k-1} + \tau (c-A\tilde{x}^{k}-Bz^{k}).
\end{align*}
Since $\psi(\lambda)=g^{\star}(B^{T}\lambda) - c^{T}\lambda$, we have $\partial \psi(\lambda)=B\partial g^{\star}(B^{T}\lambda) - c$, which implies:
\begin{align*}
\lambda^{k-1} -\tau A\tilde{x}^{k} \in \;\tau \partial \psi (\lambda^{k-1} + \tau (c-A\tilde{x}^{k}-Bz^{k})) +\lambda^{k-1} + \tau (c-A\tilde{x}^{k}-Bz^{k}).
\end{align*}
Since $z^{k}=\tilde{z}^{k}-\theta^{k}$, it follows that:
\begin{align*}
\lambda^{k-1} -\tau A\tilde{x}^{k} \in \;\tau \partial \psi (\lambda^{k-1} + \tau (c-A\tilde{x}^{k}-B\tilde{z}^{k}+B\theta^{k}))  +\lambda^{k-1} + \tau (c-A\tilde{x}^{k}-B\tilde{z}^{k}+B\theta^{k}).
\end{align*}
By Step~3 in Algorithm~\ref{al:Inexact AMA}, the above equation results in:
\begin{equation*}
\lambda^{k-1} -\tau A\tilde{x}^{k} \in \tau \partial \psi (\lambda^{k}+\tau B\theta^{k})  + \lambda^{k} +\tau B\theta^{k}.
\end{equation*}
From Step~1 in Algorithm~\ref{al:Inexact AMA} and the property of the conjugate function $p\in \partial f(q) \Leftrightarrow q\in \partial f^{\star}(p)$, we obtain:
%
\begin{equation*}
\lambda^{k-1} -\tau A(\nabla f^{\star}(A^{T}\lambda^{k})+\delta^{k}) \in \tau \partial \psi (\lambda^{k}+\tau B\theta^{k})  +\lambda^{k}+\tau B\theta^{k}.
\end{equation*}
By definition of the function $\phi$, we get:
\begin{equation*}
\lambda^{k-1} -\tau(\nabla \phi(\lambda^{k-1})+A\delta^{k}) \in \tau \partial \psi (\lambda^{k}+\tau B\theta^{k})  +\lambda^{k}+\tau B\theta^{k},
\end{equation*}
which is equivalent to:
\begin{equation*}
\lambda^{k} = \mbox{prox}_{\tau \psi}(\lambda^{k-1} - \tau (\nabla \phi(\lambda^{k-1}) + e^{k})) - \tau B\theta^{k},
\end{equation*}
with $e^{k}=A\delta^{k}$. In order to complete the proof of equation (\ref{eq: equivalence of inexact AMA and inexact ISTA on dual}), we need to show that $\lambda^{k}$ is an inexact solution of the proximal operator as defined in equation (\ref{eq:epsilon error in proximity operator}) with the error $\epsilon^{k}=\tau ^{2}L(\psi)\|B\theta^{k}\| + \frac{\tau^{2}}{2}\|B\theta^{k}\|^{2}$, i.e. to prove:
\begin{equation*}
\tau \psi(\lambda^{k}) + \frac{1}{2}\|\lambda^{k}-v\|^{2} \; \leq\; \epsilon^{k} + \mbox{min}_{\lambda} \left\lbrace \tau \psi(\lambda) + \frac{1}{2}\|\lambda-v\|^{2}\right\rbrace,
\end{equation*}
where $\nu = \lambda^{k-1} - \tau (\nabla \phi(\lambda^{k-1}) + A\delta^{k})$. Finally, using
\begin{align*}
\tau\psi(\lambda^{k}+\tau B\theta^{k})+\frac{1}{2}\|\lambda^{k}+\tau B\theta^{k}-\nu\|^{2}-\tau\psi(\lambda^{k})-\frac{1}{2}\|\lambda^{k}-\nu\|^{2} &\leq  \tau (\psi(\lambda^{k}+\tau B\theta^{k})-\psi(\lambda^{k})) + \frac{1}{2}\|\tau B\theta^{k}\|^{2} \\
&\leq  \tau ^{2}L(\psi)\|B\theta^{k}\| + \frac{\tau^{2}}{2}\|B\theta^{k}\|^{2}= \epsilon^{k},
\end{align*}
equation (\ref{eq: equivalence of inexact AMA and inexact ISTA on dual}) is proved.
\end{IEEEproof}

\subsection{Proof of Lemma~\ref{le: convergence for series for sufficient condition for IAMA with linear rate}}\label{SEC: appendix proof of le: convergence for series for sufficient condition for IAMA with linear rate}

\begin{IEEEproof}
We first prove that there exists an upper bound on the series $b^{k} = \sum^{k}_{p=1} \frac{\alpha^{-p}}{p}$. Since $0 < \alpha < 1$, there always exists a positive integer $k^{\prime}$ such that $0 < \frac{\alpha^{-k^{\prime}}}{k^{\prime}} < \frac{\alpha^{-(k^{\prime} +1)}}{k^{\prime}+1}$. We can write the series $b^{k}$ as
\begin{equation*}
b^{k} = \sum^{k^\prime}_{p=1} \frac{\alpha^{-p}}{p} + \sum^{k}_{p=k^\prime} \frac{\alpha^{-p}}{p}\enspace .
\end{equation*}
Since $k^{\prime}$ satisfies $0 < \frac{\alpha^{-k^{\prime}}}{k^{\prime}} < \frac{\alpha^{-(k^{\prime} +1)}}{k^{\prime}+1}$ and $0 < \alpha <1$, then we know that for any $t \geq k^{\prime}$ the function $\frac{\alpha^{-t }}{t}$ is a non-decreasing function with respect to $t$. Due to the fact that for any non-decreasing function $f(t)$, the following inequality holds. 
\begin{equation*}
\sum_{p\in \mathbb{Z}: y \leq p\leq x} f(p) = \int^{x}_{y} f(\lfloor t \rfloor)dt + f(x) \leq  \int^{x}_{y} f(t)dt + f(x)
\end{equation*}
where $\lfloor \cdot \rfloor$ denotes the floor operator, the series $b^{k}$ can be upper-bounded by 
\begin{equation*}
b^{k} \leq \sum^{k^\prime}_{p=1} \frac{\alpha^{-p}}{p} + \int^{k}_{k^\prime} \frac{\alpha^{-t}}{t} dt +  \frac{\alpha^{-k}}{k}\enspace .
\end{equation*}
\color{black}
We know that the integral of the function $\frac{\alpha^{-t}}{t}$ is equal to $\mathbf{E_i} (-x \log (\alpha))$,
%
%
where $ \mathbf{E_i} (\cdot)$ denotes the Exponential Integral Function, defined as $\mathbf{E_i} (x):= \int^{
\infty}_{-x} \frac{e^{-t}}{t} dt$. By using the fact that $-\mathbf{E_i} (-x) = \mathbf{E_1} (x)$, where $\mathbf{E_1} (x):= \int^{
\infty}_{x} \frac{e^{-t}}{t} dt$, and inequality~(5.1.20) in \cite{Abramowitz_Handbook_1964}, it follows that the Exponential Integral Function $\mathbf{E_i} (x)$ satisfies 
\color{black}
\begin{equation*}
-\log(1+ \frac{1}{x}) < e^{x}\cdot\mathbf{E_i} (-x) < -\frac{1}{2}\cdot \log (1+ \frac{2}{x})\enspace .
\end{equation*}
Since $e^{x} >0$, we can rewrite the inequality as
\begin{equation*}
-e^{-x}\log(1+ \frac{1}{x}) < \mathbf{E_i} (-x) < -\frac{1}{2}e^{-x} \log (1+ \frac{2}{x})\enspace .
\end{equation*}
Hence, the series $b^{k}$ can be further upper-bounded by 
\begin{align*}
b^{k} & \leq \sum^{k^\prime}_{p=1} \frac{\alpha^{-p}}{p} +  \frac{\alpha^{-k}}{k} + \mathbf{E_i} (-k \log (\alpha))- \mathbf{E_i} (-k^{\prime} \log (\alpha))\\
& < \sum^{k^\prime}_{p=1} \frac{\alpha^{-p}}{p} +  \frac{\alpha^{-k}}{k} - \frac{1}{2}e^{-k \log (\alpha)}\log (1+ \frac{2}{k \log (\alpha)})  + e^{-k^{\prime} \log (\alpha)}\log(1+ \frac{1}{k^{\prime} \log (\alpha)})\\
& = \sum^{k^\prime}_{p=1} \frac{\alpha^{-p}}{p} +  \frac{\alpha^{-k}}{k} - \frac{1}{2}\alpha^{-k}\log (1+ \frac{2}{k \log (\alpha)}) + \alpha^{-k^{\prime}}\log(1+ \frac{1}{k^{\prime} \log (\alpha)})\enspace .
\end{align*}
We can now find the upper-bound for the series $S^{k}$ as
\begin{align*}
s^{k} = \alpha^{k}\cdot b^{k} & < \alpha^{k}\sum^{k^\prime}_{p=1} \frac{\alpha^{-p}}{p} +  \frac{1}{k} - \frac{1}{2}\log (1+ \frac{2}{k \log (\alpha)}) + \alpha^{k-k^{\prime}}\log(1+ \frac{1}{k^{\prime} \log (\alpha)})\enspace .
\end{align*}
Since $0 < \alpha < 1$ and the integer $k^{\prime}$ is a constant for a given $\alpha$, the upper bound above converges to zero, as $k$ goes to infinity. In addition, we know that the two terms $\alpha^{k}\sum^{k^\prime}_{p=1} \frac{\alpha^{-p}}{p}$ and $\alpha^{k-k^{\prime}}\log(1+ \frac{1}{k^{\prime} \log (\alpha)})$ converge to zero linearly with the constant $\alpha$. From Taylor series expansion, we know that the term $\frac{1}{2}\log (1+ \frac{2}{k \log (\alpha)})$ converges to zeros at the rate $O(\frac{1}{k})$. Note that since $0 < \alpha < 1$, the term $\frac{1}{2}\log (1+ \frac{2}{k \log (\alpha)})$ is always negative for all $k>0$. To summarize, we know that the upper bound above converges to zero with the rate $O(\frac{1}{k})$. Therefore, we conclude that the series $s^{k}$ converges to zero, as $k$ goes to infinity. In addition, the convergence rate is $O(\frac{1}{k})$.
\end{IEEEproof}


%
%
%
%


%
%
%
%


\ifCLASSOPTIONcaptionsoff
  \newpage
\fi



%

\bibliographystyle{plain}
\bibliography{IFAMA_DMPC}

%

\begin{IEEEbiography}[{\includegraphics[width=1in,height=1.25in,clip,keepaspectratio]{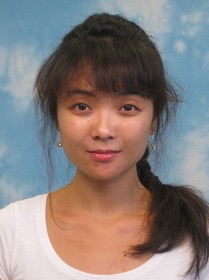}}]{Ye Pu}
received the B.S. degree from the School of Electronic Information and Electrical Engineering at Shanghai Jiao Tong University, China, in 2008, and the M.S. degree from the department of Electrical Engineering and Computer Sciences at the Technical University Berlin, Germany, in 2011. Since February 2012, she has been a Ph.D. student in the Automatic Control Laboratory at Ecole Polytechnique F\'ed\'erale de Lausanne (EPFL), Switzerland.

Her research interests are in the area of fast and distributed predictive control and optimization and distributed algorithms with communication limitations. 
\end{IEEEbiography}

\begin{IEEEbiography}[{\includegraphics[width=1in,height=1.25in,clip,keepaspectratio]{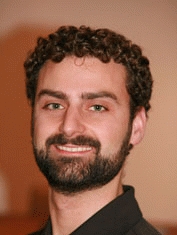}}]{Colin N. Jones}
 received the Bachelor’s degree in Electrical Engineering and the Master’s degree in Mathematics from the University of British Columbia, Vancouver, BC, Canada, and the Ph.D. degree from the University of Cambridge, Cambridge, U.K., in 2005. He is an Assistant Professor in the Automatic Control Laboratory at the École Polytechnique Fédérale de Lausanne (EPFL), Lausanne, Switzerland. He was a Senior Researcher at the Automatic Control Laboratory of the Swiss Federal Institute of Technology Zurich until 2010. 
 
 His current research interests are in the areas of high-speed predictive control and optimisation, as well as green energy generation, distribution and management.
\end{IEEEbiography}

\begin{IEEEbiography}[{\includegraphics[width=1in,height=1.25in,clip,keepaspectratio]{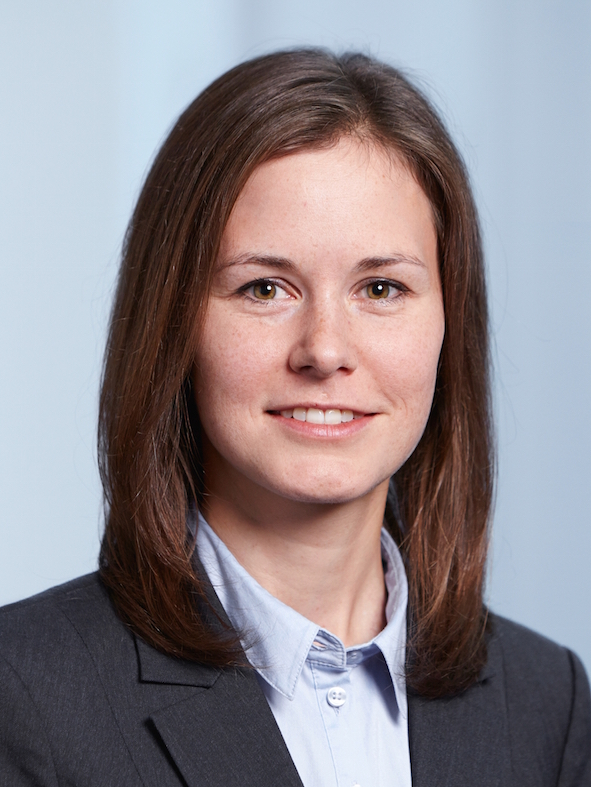}}]{Melanie N. Zeilinger}
received the Diploma degree in engineering cybernetics from the University of Stuttgart, Germany, in 2006, and the Ph.D.  degree (with honors) in electrical engineering from ETH Zurich, Switzerland, in 2011. She is an Assistant Professor at the Department of Mechanical and Process Engineering at ETH Zurich, Switzerland. She was a Marie Curie fellow and Postdoctoral Researcher with the Max Planck Institute for Intelligent Systems, Tübingen, Germany until 2015 and with the Department of Electrical Engineering and Computer Sciences at the University of California at Berkeley, CA, USA, from 2012 to 2014. From 2011 to 2012 she was a Postdoctoral Fellow with the École Polytechnique Fédérale de Lausanne (EPFL), Switzerland. 

Her current research interests include distributed control and optimization, as well as safe learning-based control, with applications to energy distribution systems and human-in-the-loop control.
\end{IEEEbiography}

%
%
%




\end{document}